\documentclass[review]{elsarticle}
\usepackage{amsthm}
\usepackage{amssymb}
\usepackage{amsmath}
\usepackage{mathrsfs}
\usepackage[colorlinks,linkcolor=blue,citecolor=magenta]{hyperref}
\usepackage{graphicx}
\newcommand{\NI}{\noindent}
\newtheorem{theorem}{\NI{\bf Theorem}}[section]
\newtheorem{lemma}{\NI\bf Lemma}[section]

\newtheorem{cor}{\NI\bf Corollary}[section]
\newtheorem{remark}{\NI\bf Remark}[section]

\newtheorem{example}{\NI\bf Example}[section]
\newcommand{\bt}{\begin{theorem}}
\newcommand{\et}{\end{theorem}}
\newcommand{\bc}{\begin{cor}}
\newcommand{\ec}{\end{cor}}
\newcommand{\bl}{\begin{lemma}}
\newcommand{\el}{\end{lemma}}
\newcommand{\br}{\begin{remark}}
\newcommand{\er}{\end{remark}}
\newcommand{\bx}{\begin{example}}
\newcommand{\ex}{\end{example}}
\newcommand{\bea}{\begin{eqnarray}}
\newcommand{\eea}{\end{eqnarray}}
\newcommand{\ben}{\begin{eqnarray*}}
\newcommand{\een}{\end{eqnarray*}}

\newcommand{\be}{\begin{eqnarray}}
\newcommand{\ee}{\end{eqnarray}}

\biboptions{authoryear}
\begin{document}
\begin{frontmatter}
\title{Equivalence of the Hazard Rate and Usual
Stochastic Orders for
Parallel Systems}


\author[mymainaddress]{Khaled Masoumifard}

\address[mymainaddress]{Department of Statistics, Faculty of Mathematical Sciences, Shahid Beheshti University, Tehran, Iran}

\begin{abstract}In this paper, we investigate stochastic comparisons of parallel
systems, and obtain two characterization results in this regard.
 First, we compare a
parallel system with independent heterogeneous components to a
parallel system with homogeneous components, and establish some
certain assumptions under which the hazard rate and usual
stochastic orders between the lifetimes of two parallel systems
are equivalent. Next, we turn our attention to two parallel
systems with their component lifetimes following multiple-outlier
model and prove that under some specified assumptions, the
$p$-larger order between the vectors of scale parameters is
equivalent to the hazard rate order as well as the usual
stochastic order between the lifetimes of these systems.
 The results established here are applicable to compute
an upper bound for the hazard rate function and a lower bound for the survival function of a parallel systems
consisting of heterogeneous components.
\end{abstract}


\end{frontmatter}

\section{Introduction}

Consider a system with $n$ components, and assume this system
fails if and only if at least $k$ components fail.  Such a system
is said to be an ($n-k+1$)-out-of-$n$ system. Series and Parallel
structures corresponding to $n$-out-of-$n$ and $1$-out-of-$n$
systems, respectively, appear extensively in the practice. In
reliability theory, there exists an intimate relation between the
($n-k+1$)-out-of-$n$ systems and the theory of order statistics.
Suppose $X_1,\cdots,X_n$ denote the lifetimes of components of a
system and $X_{1:n}\le\cdots\le X_{n:n}$ represent the
corresponding order statistics. Then,  the lifetime of an
($n-k+1$)-out-of-$n$ system comprising these components coincides
with the $k$-th order statistic $X_{k:n}$. By using the theory of
order statistics, many properties of ($n-k+1$)-out-of-$n$ systems
have been established in the literature. Order statistics also
play an importance role in various fields of probability and
statistics; see, for example,\citep{Arnold_1992} and \citep{David_Nagaraja_2003}. \\
\subsection{Scale model}
A set of independent random variables $X_1,\cdots,X_n$ with
respective distribution functions $F_1,\cdots,F_n$ is said to
follow the scale family model, if there exist positive constants
$\lambda_1,\cdots,\lambda_n$ and an absolutely continuous
distribution function $F$, with corresponding density function
$f$, such that $F_{i}(x)=F(\lambda_ix)$ for $i=1,\cdots,n$. In
this case, $F$ is called the baseline distribution, $f$ is called
the baseline density and $\lambda_1,\cdots,\lambda_n$ are scale
parameters. Let $f_i$ denote the density functions of $X_i$,
$i=1,\cdots,n$, then we have $f_i(x)=\lambda_i\,f(\lambda_ix)$,
$i=1,\cdots,n$. Some well-known distributions such as normal,
Weibull, gamma and generalized gamma are special cases of the
scale model. For this reason, the scale model has received
considerable attention in the statistical literature. One may
refer to Chapters 7 and 16 of \citep{Marshall_Olkin_2007} for more
details on the scale model and its applications.\\
\subsection{Exponentiated
scale model}
Suppose  $F(.)$ is a specified distribution with  survival
$\bar{F}(.)$, hazard rate $r(.)$ and reversed hazard rate
$\tilde{r}(.)$, and $\alpha$ is a positive number.
Set
$G(x)=[F(x)]^\alpha$. Then, it is easy to observe that $G(.)$ is
also a distribution function. In this case, we say that $F$ is
exponentiated by $\alpha$. If  $\tilde{r}_G(.)$ denotes the
reversed hazard rate function corresponding to $G$, then we have
$\tilde{r}_{G}(x)=\alpha\tilde{r}_{F}(x)$, which is referred to as
the proportional reversed hazard rate or resilience model in the
literature (see \citep{Marshall_Olkin_2007}, p. 234);

By using the exponentiation method, one can  produce a new
distribution with some attractive properties. For example, suppose
$F(x)=1-e^{-x}$, $x>0$, and $\alpha$ is a positive number. Then,
$G(x)=(1-e^{-x})^\alpha$, $x>0$, is also a distribution function,
which is called the generalized exponential distribution. It is
well-known that the exponential distribution has constant hazard
rate whereas the generalized exponential distribution has
decreasing hazard rate for $\alpha\le 1$ and increasing hazard
rate for $\alpha\ge 1$. For a comprehensive discussion on the
generalized exponential distribution, one may refer to \citep{Gupta_Kundu_(1999)}, \citep{Gupta_Kundu_(2004)}.

Suppose  $F(.)$ is an absolutely continuous distribution function
with corresponding density $f(.)$. Then, independent random
variables $X_1,\cdots,X_n$ are said to belong to the
exponentiated scale family of distributions if $X_i\sim
\Big(F(\lambda_i\,x)\Big)^{\alpha_i}$, where $\alpha_i>0$ and
$\lambda_i>0$ for $i=1,\cdots,n$. In this case, $F(.)$ is said to
be the baseline distribution function and  $\alpha_i$'s and
$\lambda_i$'s are the shape and  scale parameters, respectively.
Note that we can present the parameters of the exponentiated
scale model in a matrix format. 

Let us quickly recall some common notions of stochastic order, and
majorization and $p$-larger orders that will be used to establish
the main results of the paper. Throughout, we use `increasing' to
state `non-decreasing' and similarly `decreasing' to state
`non-increasing'. Suppose $X$ and $Y$ are two non-negative random
variables with distribution functions $F$ and $G,$ survival
functions $\bar{F}$ and $\bar{G},$ density functions $f$ and $g,$
and hazard rates $r_{X}=f/\bar{F}$ and $r_Y=g/\bar{G},$
respectively. Then, $X$ is said to be larger than $Y$ in the
usual stochastic order (shown by $X\ge_{st}Y$) iff
$\bar{F}(x)\ge\bar{G}(x)$ for all $x\in\mathbb{R}^+$. $X$ is said
to be larger than $Y$ in the hazard rate order (shown by $X\ge
_{hr}Y$) iff $\bar{F}(x)/\bar{G}(x)$ is increasing in
$x\in\mathbb{R}^+$. In fact, $X\ge_{hr}Y$ iff $r_{Y}(x)\ge
r_{X}(x)$ for all $x\in\mathbb{R}^+$. It is well-known that the
hazard rate order implies the usual stochastic order. Interested
readers may refer to \citep{Muller_Stoyan_(2002)} and \citep{Shaked_Shanthikumar_(2007)} for comprehensive discussions on stochastic
orders and their applications.

For two vectors $\mbox{\boldmath$x$}=(x_1,\cdots,x_n)$ and
$\mbox{\boldmath$y$}=(y_1,\cdots,y_n)$, let $\{x_{(1)},\cdots,
x_{(n)}\}$ and $\{y_{(1)},\cdots, y_{(n)}\}$ denote the
increasing arrangements of their components, respectively. Then,
the vector $\mbox{\boldmath$x$}$ is said to majorize another
vector $\mbox{\boldmath$y$}$ (written
$\mbox{\boldmath$x$}\stackrel{m}{\succeq}\mbox{\boldmath$y$}$) if
$\,\sum_{j=1}^{i}\,x_{(j)}\le\sum_{j=1}^{i}\,y_{(j)}$ for
$i=1,\cdots,n-1$, and
$\sum_{j=1}^{n}\,x_{(j)}=\sum_{j=1}^{n}\,y_{(j)}$. The vector
$\mbox{\boldmath$x$}$ in $\mathbb{R^+}^n$ is said to be $p$-larger
than another vector  $\mbox{\boldmath$y$}$ in $\mathbb{R^+}^n$
(written
$\mbox{\boldmath$x$}\stackrel{p}{\succeq}\mbox{\boldmath$y$}$) if
$\,\prod_{j=1}^{i}\,x_{(j)}\le\prod_{j=1}^{i}\,y_{(j)}$  for
$i=1,\cdots,n$. When $\mbox{\boldmath$x$}$,$\mbox{\boldmath$y$}$
in $\mathbb{R^+}^n$,
$\mbox{\boldmath$x$}\stackrel{m}{\succeq}\mbox{\boldmath$y$}$
implies
$\mbox{\boldmath$x$}\stackrel{p}{\succeq}\mbox{\boldmath$y$}$. The
converse is, however, not true. For example,
$(1,5.5)\stackrel{p}{\ge}(2,3)$, but clearly the majorization
order does not hold between these two vectors. A real-valued
function $\phi$, defined on a set
$\mathbb{A}\subseteq\mathbb{R}^n$, is said to be Schur-convex
(Schur-concave) on $\mathbb{A}$ if
$\mbox{\boldmath$x$}\stackrel{m}{\succeq}\mbox{\boldmath$y$}$
implies
$\phi(\mbox{\boldmath$x$})\ge(\le)\phi(\mbox{\boldmath$y$})$ for
any $\mbox{\boldmath$x$},\mbox{\boldmath$y$}\in\mathbb{A}$. For
comprehensive discussions about the $p$-larger and majorization
orders with their applications, we refer the readers to Khaledi
and \citep{Kochar_(2002)} and \citep{Marshall_et_al_(2011)}.


 In this paper, we compare parallel systems with respect to the
hazard rate and usual stochastic orders when their components
follow the exponentiated scale family of distributions. First, suppose $X_1,\cdots,X_n$ are
independent non-negative random variables with $X_i\sim
\bigg(F(\lambda_ix)\bigg)^{\alpha_i}$, $i=1,\cdots,n$, and suppose $Y_1,\cdots,Y_n$ are
independent non-negative random variables with common
distribution function $\bigg(F(\lambda x)\bigg)^{\alpha_i}$. Then, it is shown under
some certain assumptions that \bea
\lambda \ge \lambda_{wg} \Longleftrightarrow
X_{n:n}\ge_{hr}Y_{n:n}\Longleftrightarrow
X_{n:n}\ge_{st}Y_{n:n},
\label{e1}\eea 
where
$\lambda_{wg}=\prod_{i=1}^{n}\lambda_i^{{\alpha_i}/{n\bar{\alpha}}}$ is the weighted
geometric mean of $\lambda_i$'s. We show that the results in
(\ref{e1}) hold for the power generalized weibull and exponentiated generalized gamma family of
distributions. By means of a counterexample, we
show that the result in (\ref{e1}) cannot be extended to the general
case when the lifetimes of two parallel systems are heterogeneous
and the vectors of scale parameters are ordered with respect to
the $p$-larger order. So, for extending (\ref{e1}), in the first step
we turn our attention to the multiple-outlier scale model which
reduces the complete heterogenity of the scale parameters.
Specifically, suppose $X_1,\cdots,X_n$ are independent
non-negative random variables with $X_i\sim F(\lambda x)$,
$i=1,\cdots,p$, and $X_j\sim F(\lambda ^*x)$, $j=p+1,\cdots,n$.
Further, suppose $Y_1,\cdots,Y_n$ is another set of independent
non-negative random variables with $Y_i\sim F(\mu x)$,
$i=1,\cdots,p$, and $Y_j\sim F(\mu^*x)$, $j=p+1,\cdots,n$. Assume
that $\lambda\le\mu\le\mu^*\le\lambda^*$. Then, under some
specified assumptions, we show that \bea
(\underbrace{\lambda,\cdots,\lambda}_{p},
\underbrace{\lambda^*,\cdots,\lambda^*}_{n-p})\stackrel{p}{\succeq}
(\underbrace{\mu,\cdots,\mu}_{p},
\underbrace{\mu^*,\cdots,\mu^*}_{n-p}) \Longleftrightarrow
X_{n:n}\ge_{hr}Y_{n:n}\Longleftrightarrow
X_{n:n}\ge_{st}Y_{n:n}\label{e2}.\eea We also show that the results in
(\ref{e2}) hold for the power generalized weibull and generalized gamma family of
distributions. The results established here reinforce and extend
the well-known results  in \citep{Zhao_(2011)},\citep{Zhao_Balakrishnan_(2012)}, Balakrishnan and \citep{Zhao_(2013a)} and \citep{Balakrishnan_et_al_(2014)} which deal with the exponential, gamma and generalized exponential
distributions.

\section{Preliminaries}
In this section, we discuss some aspects of the scale model and
present some useful lemmas that will be used in the sequel.
Suppose $X_1,\cdots,X_n$ are independent non-negative random
variables following the scale model with the scale parameters
$\lambda_1\cdots,\lambda_n,$ the baseline distribution $F,$ and
the baseline density $f$.
 In this paper, we assume some
certain assumptions on the baseline density to simplify the
computation. More precisely, suppose the baseline density function
can be rewritten as \bea f(x)=w(x)\,h(x),\qquad\quad
x\in\mathbb{R^+},\label{e3}\eea where $w(x)$ and $h(x)$ are two
differentiable
positive functions with the following properties:\\
{\bf C-1.}\,\,$w(xy)=w(x)\,w(y)$ for all $x,y\in\mathbb{R}^+$;\\
{\bf C-2.}\,\,$\lim_{x\rightarrow 0} \frac{\displaystyle h(\lambda_i x)}{\displaystyle h(\lambda_j x)}=1$.\\ As an example of the decomposition form in
(\ref{e3}),
\begin{itemize}
\item
set
 \ben w(x)=x^{p-1}\qquad and\qquad
h(x)=\frac{\displaystyle p}{\displaystyle q}(1+x^p)^{\frac{1}{q}-1}e^{1-(1+x^p)^{\frac{1}{q}}},\qquad
x\in\mathbb{R^+},\een
 where $p>0$, $q>0$. Clearly,
the assumptions in {\bf C-1} and {\bf C-2} satisfy for this case.
Under this setting, the baseline density function reduces to the
density function of the power generalized weibull distribution with shape
parameters $p$ and $q$, and scale parameter 1 (on in
short $PGW(p,q,1)$). The power generalized weibull distribution has a decreasing failure rate when $p\le q,\,p\le 1$, an increasing failure rate when $p\ge q,\,p\ge 1$, a bathtub failure rate when
$0<q<p<1$ and an upside down bathtub (or unimodal) failure rate when $q>p>1$.
It includes Weibull and exponential 
distributions as special cases. For more details on this family and its applications in probability and statistics,the reader is
referred to \citep{Bagdonavicius} and \citep{Nikulin_(2002)}.
\item
 set \ben w(x)=x^{\beta-1}\qquad and\qquad
h(x)=\frac{\displaystyle\alpha}{\displaystyle\Gamma^*(\frac
{\displaystyle\beta}{\displaystyle\alpha})}\,e^{-x^{\alpha}},\qquad
x\in\mathbb{R^+},\een where $\alpha>0$, $\beta>0$, and
$\Gamma^*(.)$ indicates the gamma function. Clearly,
the assumptions in {\bf C-1} and {\bf C-2} satisfy for this case.
Under this setting, the baseline density function reduces to the
density function of the generalized gamma distribution with shape
parameters $\alpha$ and $\beta$, and scale parameter 1 (on in
short $GG(\alpha,\beta,1)$). The generalized gamma distribution
includes exponential, Weibull and gamma distributions  as special
cases. Moreover, the log-normal distribution can be also obtained
from the generalized gamma distribution as $\alpha=2$ and
$\beta\rightarrow\infty$. The generalized gamma distribution also
has a flexible hazard rate function which is  increasing for
$\alpha\ge 1$ and $\beta\ge 1$, decreasing for $\alpha\le 1$ and
$\beta\le 1$, bathtub shape for $\alpha<1$ and $\beta>1$ and
upside-down bathtub shape for $\alpha>1$ and $\beta<1$. For more
details on the generalized gamma distribution and its
applications, we refer the readers to \citep{Balakrishnan_and_Peng_(2006)} and \citep{Kleiber_Kotz_(2003)}.
\end{itemize}
From the assumptions in {\bf C-1} and {\bf C-2}, we easily obtain
the following result. \br
$w(1)=1$ and $w(x)=x^{w'(1)}$ for every
$x\in\mathbb{R}^+$, where $w'(x)$ denote the derivative of $w(x)$
with respect to $x$.
\er
In the following lemma, we compare a parallel system consisting of
components whose lifetimes follow the exponentiated scale model to the same
system with lifetimes of its components following the homogeneous exponentiated
scale model, and establish some certain assumptions under which
the hazard rate and usual stochastic orders between the lifetimes
of these systems are equivalent.
\begin{lemma}
Suppose $X_1,\cdots,X_n$ are independent non-negative random
variables with $X_i\sim \big(F(\lambda_i\,x)\big)^{{\alpha}_i}$, $i=1\cdots,n$. Further,
suppose $Y_1,\cdots,Y_n$ are independent non-negative random
variables with  $Y_i\sim \big(F(\lambda x)\big)^{{\alpha}_i}$.
Assume that following assumptions hold:
\begin{itemize}
\item[($a$)]$x\tilde{r}(x)$ is decreasing in $x\in\mathbb{R}^+$;
\item[($b$)]$x(\tilde{r}(x)-\frac{h'(x)}{h(x)})$ is increasing
in $x\in\mathbb{R}^+$;
\item[($c$)]$\bar{F}(x)\le -\frac{w(x)\,h^{2}(x)}{h'(x)}$ for every
$x\in\mathbb{R}^+$.
\end{itemize}
Where $\tilde{r}(x)={f(x)}/{F(x)}$. Then, for $\alpha_i \ge 1,\,i=1,\cdots,n$ the following three statements are equivalent:
\begin{itemize}
\item[(i)]$\lambda\ge\lambda_{wg}$;
\item[(ii)]$X_{n:n}\ge_{hr}Y_{n:n}$;
\item[(ii)]$X_{n:n}\ge_{st}Y_{n:n}$.
\end{itemize}
Where ${\lambda}_{wg}=\prod_{i=1}^n \lambda_i^{\frac{\alpha_i}{\sum_{i=1}^n \alpha_i}}$.
\label{L_new}
\end{lemma}
 Now, we derive some ordering results between parallel
systems with the lifetimes of their components following the exponentiated
scale model under the decomposition form in (\ref{e3}) and the
assumptions in {\bf C-1} and {\bf C-2}.  \\
 \begin{lemma}
Suppose $X_1,\cdots,X_n$ are independent non-negative random
variables with $X_i\sim F(\lambda x)$, $i=1,\cdots,p$, and
$X_j\sim F(\lambda^*x)$, $j=p+1,\cdots,n$. Further, suppose
$Y_1,\cdots,Y_n$ is another set of independent non-negative
random variables with $Y_i\sim F(\mu x)$, $i=1,\cdots,p$, and
$Y_j\sim F(\mu^*x)$, $j=p+1,\cdots,n$. Assume that following
assumptions hold:
\begin{itemize}
\item[($a$)]$x\tilde {r}(x)$ is decreasing in
$x\in\mathbb{R}^+$;
\item[($b$)]$x(\tilde{r}(x)-\frac{ h'(x)}{h(x)})$ is
increasing in $x\in\mathbb{R}^+$;
\item[($c$)]$\bar{F}(x)\le -\frac{w(x)\,h^2(x)}{h'(x)}$
for every $x\in\mathbb{R}^+$.
\end{itemize}
So, if $\lambda\le\mu\le\mu^*\le\lambda^*$ and
$\lambda^p{\lambda^*}^{n-p}=\mu^p{\mu^*}^{n-p}$, then
$X_{n:n}\ge_{hr}Y_{n:n}.$
\label{l-new2}
\end{lemma}
\begin{lemma}
Suppose $X_1,\cdots,X_n$ are independent non-negative random
variables with $X_i\sim F(\lambda x)$, $i=1,\cdots,p$, and
$X_j\sim F(\lambda^*x)$, $j=p+1,\cdots,n$. Further, suppose
$Y_1,\cdots,Y_n$ is another set of independent non-negative
random variables with $Y_i\sim F(\mu x)$, $i=1,\cdots,p$, and
$Y_j\sim F(\mu^*x)$, $j=p+1,\cdots,n$. Assume that
$\lambda\le\mu\le\mu^*\le\lambda^*$ and following assumptions
hold:
\begin{itemize}
\item[($a$)]$x\tilde {r}(x)$ is decreasing in
$x\in\mathbb{R}^+$;
\item[($b$)]$x(\tilde{r}(x)-\frac{ h'(x)}{h(x)})$ is
increasing in $x\in\mathbb{R}^+$;
\item[($c$)]$\bar{F}(x)\le -\frac{w(x)\,h^2(x)}{h'(x)}$
for every $x\in\mathbb{R}^+$.
\end{itemize}
Then, the following three statements are equivalent:
\begin{itemize}
\item[(i)]$(\underbrace{\lambda,\cdots,\lambda}_{p},
\underbrace{\lambda^*,\cdots,\lambda^*}_{n-p})\stackrel{p}{\succeq}
(\underbrace{\mu,\cdots,\mu}_{p},
\underbrace{\mu^*,\cdots,\mu^*}_{n-p})$;
\item[(ii)]$X_{n:n}\ge_{hr}Y_{n:n}$;
\item[(iii)]$X_{n:n}\ge_{st}Y_{n:n}$.
\end{itemize}
\label{l-new3}
\end{lemma}

\section{Comparison of $X_{n:n}$ from Heterogenous and Homogeneous Samples}
Stochastic comparisons of a
parallel system with independent heterogeneous components to a 
parallel system with homogeneous components were first made by \citep{Dykstra_et_al}
Let $(X_1,\cdots,X_n)$ and $(Y_1,\cdots,Y_n)$ be two independent random vectors,
\begin{itemize}
\item If $X_i\sim E(\lambda_i)$\footnote{Exponential distribution with hazard rate $\lambda_i$} and $Y_i\sim E({\lambda})$, $i=1,\cdots, n$, \citep{Dykstra_et_al} showed that
\begin{eqnarray}
\lambda=\bar{\lambda}\Longrightarrow X_{n:n}\ge_{hr} Y_{n:n}
\label{dykstra}
\end{eqnarray}
where $\bar{\lambda}=\frac{1}{n}\sum_{i=1}^n \lambda_i$, the arithmetic mean of $\lambda_i$'s.  Which was further strengthened by \citep{kochar2007some} as
$$
 X_{n:n}\ge_{lr} Y_{n:n}.
 $$
 [10] also strengthened the result in (\ref{dykstra}), under a weaker condition,
 $$
\lambda=\tilde{\lambda}\Longrightarrow X_{n:n}\ge_{hr} Y_{n:n}
$$
where $\tilde{\lambda}=\bigg(\prod_{i=1}^n \lambda_i\bigg)^{1/n}$, the geometric mean of $\lambda_i$'s.
\item If $X_i\sim G(r,\lambda_i)$\footnote{Gamma distribution with shape
parameter $r$ and scale parameter $\lambda_i$} and $Y_i\sim G(r,\tilde{\lambda})$, $i=1,\cdots, n$, \citep{balakrishnan2013hazard} show that, for $0<r\le 1$,
$$
 X_{n:n}\ge_{hr} Y_{n:n}.
$$
  
\item If $X_i\sim GE(\alpha_i,\lambda_i)$\footnote{Generalized Exponential distribution with shape
parameter $\alpha_i$ and scale parameter $\lambda_i$} and $Y_i\sim GE(\bar{\alpha},{\lambda_{wg}})$, $i=1,\cdots, n$,  \citep{balakrishnan2015stochastic} show that, for $\bar{\alpha}\ge 1$,
 $$
 X_{n:n}\ge_{hr} Y_{n:n}
$$ 
 where $\bar{\alpha}=\frac{1}{n}\sum_{i=1}^n \alpha_i$ and $\lambda_{wg}=\prod_{i=1}\lambda_i^{{\alpha_i}/{n\bar{\alpha}}}$. They also showed that, If $X_i\sim GE(\alpha,\lambda_i)$ and $Y_i\sim GE({\alpha},\bar{\lambda})$, $i=1,\cdots, n$,  then
   $$
 X_{n:n}\ge_{lr} Y_{n:n}.
 $$
   \end{itemize}

The case of power-generalized weibull distribution and exponentiated generalized gamma distribution are discussed in the following subsections, respectively.
\subsection{Power-generalized weibull distribution}
In the following lemma, we show that the assumptions of Lemma \ref{L_new}
hold for the power-generalized weibull distribution  under a restriction
on its shape parameter.
\bl  
\label{lemma_3.1}
 Suppose the baseline distribution is
$PGW(p,q,1)$. Then, for $q\le 1$,  the
assumptions of Lemma \ref{L_new} hold.
\label{L-pwd}
\el
 
 Now, from  Lemma \ref{L_new} and Lemma \ref{L-pwd}, we obtain the following Theorem.
\begin{theorem} 
 Let $(X_{1},\cdots,X_{n})$ and $(Y_{1},\cdots,Y_{n})$ be two independent random vectors with
$X_{i} \sim PGW(p,q, \lambda_i)$ and  $Y_{i} \sim PGW(p,q,\bar{\lambda})$,  $i=1,\cdots, n$.
Where $\bar{\lambda}=\big(\prod_{i=1}^n \lambda_i\big)^{{1}/{n}}$.
Then, for $q \le 1 $, the following three statements are
equivalent:
\begin{itemize}
\item[(i)]$\lambda\ge \bar{\lambda}$;
\item[(ii)]$X_{n:n}\ge_{hr}Y_{n:n}$;
\item[(ii)]$X_{n:n}\ge_{st}Y_{n:n}$.
\end{itemize}
\label{T-pgw}
\end{theorem}
 
\subsection{Exponentiated generalized gamma distribution}
\cite{cordeiro2011exponentiated} then employed the exponentiation method on
the generalized  gamma  distribution to introduce a four-parameter
lifetime distribution, called the exponentiated generalized gamma
distribution. Thus, a random variable $X$ is said to have the
exponentiated generalized gamma distribution with shape
parameters ${\gamma}$, $\alpha$ and $\beta$, and scale parameter
$\lambda$ (denote by $X\sim EGG({\gamma},\alpha,\beta,\lambda)$) if
its cumulative distribution function is given by \ben
F(t;{\gamma},\alpha,\beta,\lambda)=
\Bigg[\gamma\Big((\lambda\,t),{\alpha},{\beta}\Big)\Bigg]^{\theta},\qquad
t>0, {\theta}>0, \alpha>0, \beta>0, \lambda>0,\een  where
$\gamma(x,{\alpha},{\beta})$ is the cumulative distribution function of a
generalized gamma distribution with shape parameters ${\alpha}$ and ${\beta}$ and scale
parameter 1.  Many well-known distributions are sub-models of the
exponentiated generalized gamma distribution. For ${\alpha}={\beta}$, it
becomes the exponentiated Weibull distribution proposed by
\citep{cordeiro2011exponentiated} If ${\alpha}=\beta=1$, it reduces
to the generalized exponential distribution introduced by \citep{cordeiro2011exponentiated} For $\gamma=1$ and  ${\alpha}=\beta=2$, it
becomes the Rayleigh distribution. When $\beta=1$,  it reduces to
exponentiated gamma distribution initiated by \citep{gupta1998modeling} If ${\alpha}/{\beta}=\gamma=1$, the two-parameter Weibull
distribution is obtained, while for ${\gamma}=\alpha=\beta=1$ the
exponential distribution is deduced. An interesting  property of
the exponentiated generalized gamma distribution is that its
hazard rate admits bathtub, upside-down bathtub or monotone
shapes. 
For more details on some general properties of the exponentiated
generalized gamma distribution and its applications, one may
refer to \citep{cordeiro2011exponentiated}\\
In the following lemma, we show that the assumptions of Lemma \ref{L_new}
hold for the exponentiated 
generalized gamma distribution under a restriction
on its shape parameter.
\bl 
 Suppose the baseline distribution is
$GG(\alpha,\beta,1)$. Then, for $\alpha\ge\beta$,  the
assumptions of Lemma \ref{L_new} hold.
\label{L-pwd1}
\el
 
The following Theorem is a direct consequence of  Lemma \ref{L_new}
and Lemma \ref{L-pwd1}.
 \begin{theorem}
 	\label{Theo_3.2}
 	 Suppose $X_1,\cdots,X_n$ are independent non-negative random
variables with $X_i\sim  EGG({\gamma}_i,\alpha,\beta,\lambda_i)$, $i=1\cdots,n$. Further,
suppose $Y_1,\cdots,Y_n$ are independent non-negative random
variables with  $Y_i\sim  EGG({\gamma}_i,\alpha,\beta,\lambda)$.
Then, for $\alpha \ge \beta$ and $\gamma_i \ge 1,\,i=1,\cdots,n$ the following three statements are equivalent:
\begin{itemize}
\item[(i)]$\lambda\ge \lambda_{wg}$;
\item[(ii)]$X_{n:n}\ge_{hr}Y_{n:n}$;
\item[(ii)]$X_{n:n}\ge_{st}Y_{n:n}$.
\end{itemize}
Where ${\lambda}_{wg}=\prod_{i=1}^n \lambda_i^{\frac{\alpha_i}{\sum_{i=1}^n \alpha_i}}$.
\label{T-EGG}
\end{theorem}
Theorem \ref{T-pgw} and Theorem \ref{T-EGG} can be used in the practical situations
to obtain an upper bound for the hazard rate function and a lower
bound for the survival function of a parallel system.
\br
\begin{enumerate}
	\item[(i)]
	Under the conditions of Theorem \ref{T-pgw},
	\begin{itemize}
		\item the hazard $r_{X_{n:n}}$ of $X_{n:n}$ satisfies
		$$ r_{X_{n:n}}(x)\le
		\frac{\displaystyle n\frac{p}{q}\bar{\lambda}^px^{p-1}(1+(\bar{\lambda} x)^p)^{1/q-1}e^{n\{1-(1+(\bar{\lambda} x)^p)^{1/q}\}}}
		{\displaystyle 1-e^{n\{1-(1+(\bar{\lambda} x)^p)^{1/q}\}}},\qquad x>0 $$
		\item $$p(X_{n:n}\ge x+t\mid X_{n:n}\ge x) \ge \frac{\displaystyle 1-e^{n\big\{1-\big(1+\big(\bar{\lambda} (x+t)\big)^p\big)^{1/q}\big\}}}{\displaystyle 1-e^{n\big\{1-\big(1+\big(\bar{\lambda} x\big)^p\big)^{1/q}\big\}}},\qquad x>0, t\ge 0$$
	\end{itemize}
	\item[(ii)]
	Under the conditions of Theorem \ref{T-EGG},
	\begin{itemize}
		\item the hazard $r_{X_{n:n}}$ of $X_{n:n}$ satisfies
		$$
		r_{X_{n:n}}(x)\le 
		\frac{\displaystyle n\bar{\theta} \gamma'\Big(\lambda_{wg}\,x,{\alpha},{\beta}\Big) \Big[\gamma\Big(\lambda_{wg}\,x,{\alpha},{\beta}\Big)\Big]^{n\bar{\theta}-1}}
		{1-\Big[\gamma\Big(\lambda_{wg}\,x,{\alpha},{\beta}\Big)\Big]^{n\bar{\theta}}},\qquad x>0 $$
		where $\gamma'\Big(\lambda_{wg}\,x,{\alpha},{\beta}\Big)=\frac{\partial}{\partial x}\gamma\Big(\lambda_{wg}\,x,{\alpha},{\beta}\Big)$ and $\gamma(x,{\alpha},{\beta})$ is the cumulative distribution function of a
		generalized gamma distribution with shape parameters ${\alpha}$ and ${\beta}$ and scale
		parameter 1.
		\item $$p(X_{n:n}\ge x+t\mid X_{n:n}\ge x) \ge \frac{\displaystyle 1-\Big[\gamma\Big(\lambda_{wg}\,(x+t),{\alpha},{\beta}\Big)\Big]^{n\bar{\theta}}}{\displaystyle 1-\Big[\gamma\Big(\lambda_{wg}\,x,{\alpha},{\beta}\Big)\Big]^{n\bar{\theta}}},\qquad x>0,\,t\ge 0$$
	\end{itemize}
\end{enumerate}
\label{bound}
\er
So, if we want to replace the components of a
parallel system consisting of $n$ independent heterogeneous components
by identical components, then the bounds in Remark \ref{bound} enable us
to determine the exact values of the shape and scale parameters of
the lifetimes of the new components, in order to preserve the reliability.
This problem can utilize in precise production constraints of components
of a parallel system. We shall now illustrate the above observations
by a numerical example.
\bx Suppose $X_1,X_2,X_3$ are independent
random variables with $X_i\sim PGW(p,q,\lambda_i)$,
$i=1,2,3,$ where $(\lambda_1,\lambda_2,\lambda_3)=(1.5,2,3.5)$.
Further, suppose $Y_1,Y_2,Y_3$ is another set of independent
random variables with $Y_i\sim PGW(p,q,\lambda)$,
$i=1,2,3.$ Assume that the hazard rates of $X_{3:3}$ and $Y_{3:3}$ are denoted, respectively,
by $h(.;1.5,2,3.5)$ and
$h(.;\lambda,\lambda,\lambda)$.  Fig.1 and Fig.2 represent
the hazard rate functions of $X_{3:3}$ and $Y_{3:3}$ for the case
$(p,q)=(1.5,0.8)$ and $(p,q)=(0.8,0.4)$,
respectively, when  parameter $\lambda$ is taken as
$\bar{\lambda}=2.333$ (the arithmetic mean of $\lambda_i$'s) and
$\tilde{\lambda}=2.189$ (the geometric means of $\lambda_i$'s).
It appears from these figures that the hazard rate function of
$Y_{3:3}$ in terms of the geometric mean of $\lambda_i$'s
provides an upper bound for the hazard rate function of $X_{3:3}$
which is sharper than those in terms of the arithmetic mean of
$\lambda_i$'s. $\qquad\square$\ex
\begin{figure}
	\label{fig:cc}
	\begin{center}
		\includegraphics[scale=0.5]{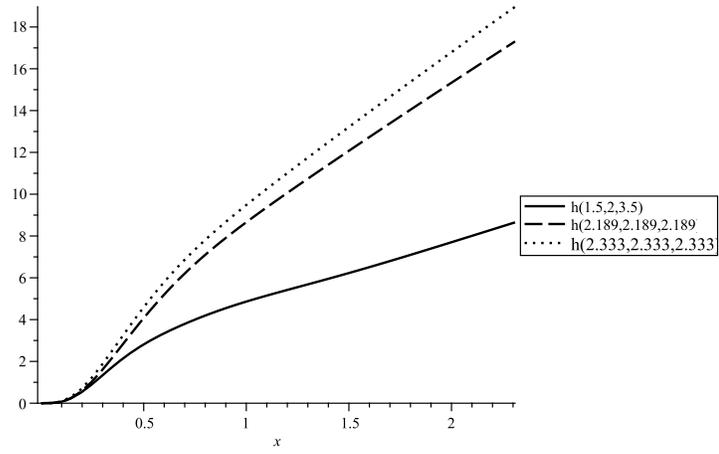}
	\end{center}
	\caption{Plot of the hazard rate functions for $p=1.5$ and $q=0.8$}
	\label{fig1}
\end{figure}
\begin{figure}
	\label{fig:cc}
	\begin{center}
		\includegraphics[scale=0.5]{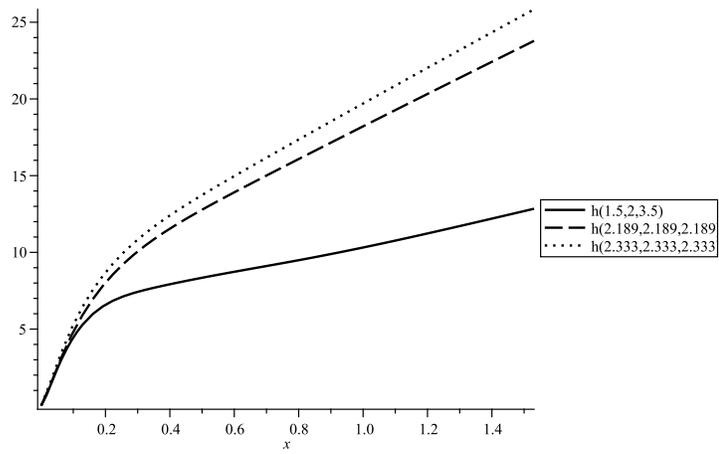}
	\end{center}
	\caption{Plot of the hazard rate functions for $p=0.8$ and $q=0.4$}
	\label{fig2}
\end{figure}

\section{Comparisons in multiple-outlier model}
It is natural to ask whether the result of Theorem \ref{T-pgw} can be
extended to the general case when the lifetimes of two parallel
systems are heterogeneous and the vectors of scale parameters are
ordered with respect to the $p$-larger order. The following
example provides a negative answer to this problem. 
\bx Suppose
$X_1,X_2,X_3$ are independent random variables with $X_i\sim
PGW(2,1,\lambda_i)$, $i=1,2,3$, where
$(\lambda_1,\lambda_2,\lambda_3)=(0.1,1,9)$. Further, suppose
$Y_1,Y_2,Y_3$ is another set of independent random variables with
$Y_i\sim PGW(2,1,\mu_i)$, $i=1,2,3$, where
$(\mu_1,\mu_2,\mu_3)=(0.1,4,6)$. It is easy to observe that
$(\lambda_1,\lambda_2,\lambda_3)\stackrel{p}{\succeq}(\mu_1,\mu_2,\mu_3)$.
On the other hand, we have
\ben\frac{\displaystyle\bar{F}_{X_{3:3}}(0.8)}{\displaystyle\bar{F}_{Y_{3:3}}(0.8)}=1.00339,
\qquad
\frac{\displaystyle\bar{F}_{X_{3:3}}(1)}{\displaystyle\bar{F}_{Y_{3:3}}(1)}=1.0037,
\qquad
\frac{\displaystyle\bar{F}_{X_{3:3}}(1.9)}{\displaystyle\bar{F}_{Y_{3:3}}(1.9)}=1.0019.\een
 Thus, the ratio
$\frac{\displaystyle\bar{F}_{X_{3:3}}(x)}{\displaystyle\bar{F}_{X_{3:3}}(x)}$
is not monotone in $x\in\mathbb{R^+}$, which means that $X_{3:3}$
and $Y_{3:3}$ can not be compared with respect to the hazard rate
order. $\qquad\square$\ex

%

Since Theorem \ref{T-pgw} could not be extended to the general case, in
the first step attention moves to the multiple-outlier scale model
which reduces the complete heterogenity of the scale parameters.
In the sequel, we consider parallel systems with the lifetimes of
their components following the multiple-outlier scale model and
investigate the problem of stochastic comparisons between these
systems. 


Let us now recall some results in the literature that are most
pertinent to the main results that are established in the sequel.
Let $Z_1,\cdots,Z_n$  be independent  exponential random variables
with  hazard rates $\lambda_1,\cdots,\lambda_n$, respectively.
Further, let $Z^{*}_1,\cdots,Z^{*}_n$ be another set of independent
exponential random variables with hazard rates
$\lambda^{*}_1,\cdots,\lambda^{*}_n$, respectively. Then, \citep{pledger1971comparisons} proved
that\bea(\lambda_1,\cdots,\lambda_n)\stackrel{m}{\succeq}
(\lambda^{*}_1,\cdots, \lambda^{*}_n) \Longrightarrow
Z_{k:n}\ge_{st} Z^*_{k:n}.\label{4.1}\eea In the special case of $k=2$ and
$n=2$, Boland et al. (1994) strengthened (\ref{4.1}) from the usual
stochastic order to the hazard rate order. Moreover, they showed, by means
 of a counterexample, that (3.10) cannot be extended
to the hazard rate order for $n\ge 3$. Recently, Zhao and
Balakrishnan (2012) improved (3.10) from the usual stochastic
order to the hazard rate order for parallel systems in
multiple-outlier exponential models. Specifically, let
$Y_1,\cdots,Y_n$ be independent exponential random variables with
$Y_1,\cdots,Y_p$ having common hazard rate $\lambda_1$ and
$Y_{p+1},\cdots,Y_n$ having common hazard rate $\lambda_2$. Further,
let $Y^*_{1},\cdots,Y^*_{n}$ be another set of independent exponential random
variables with $Y^*_{1},\cdots,Y^*_{p}$ having common hazard rate
$\lambda^{*}_{1}$ and $Y^*_{p+1},\cdots,Y^*_{n}$ having common
hazard rate $\lambda^{*}_{2}$. Then, for
$\lambda_{1}\le\lambda^{*}_{1} \le\lambda^{*}_{2}\le\lambda_{2}$,
Zhao and Balakrishnan (2012) showed that
\bea(\underbrace{\lambda_1,\cdots,\lambda_1}_{p},
\underbrace{\lambda_2,\cdots,\lambda_2}_{n-p})\stackrel{p}{\succeq}
(\underbrace{\lambda^{*}_1,\cdots,\lambda^{*}_1}_{p},
\underbrace{\lambda^{*}_2,\cdots,\lambda^{*}_2}_{n-p})\Longleftrightarrow
Y_{n:n}\ge_{hr}Y^{*}_{n:n}\Longleftrightarrow
Y_{n:n}\ge_{st}Y^{*}_{n:n}.\eea
Balakrishnan and Zhao
(2013a) and Balakrishnan et al. (2014) extended  a
result similar to the one in (3.11) for multiple-outlier gamma model and multiple-outlier
GE model, respectively.
  
From Lemma \ref{l-new3}, Lemma \ref{L-pwd} and Lemma \ref{L-pwd1} , we obtain a
result similar to the one in (3.11) for multiple-outlier
PGW model and multiple-outlier
GG model.   \bt   Let $(X_{1,1},\cdots,X_{1,n_1},X_{2,1},\cdots,X_{2,n_2})$ and $(Y_{1,1},\cdots,Y_{1,n_1},Y_{2,1},\cdots,Y_{2,n_2})$ be two independent random vectors with
 \begin{itemize}
 \item[(i)]$X_{i,j} \sim PGW(p,q, \lambda_i)$ and  $Y_{i,j} \sim PGW(p,q,{\lambda}_{i}^{*})$, $j=1,\cdots,n_i$, $i=1,2$ or
 \item[(ii)]$X_{i,j} \sim GG(\alpha,q \alpha, \lambda_i)$ and  $Y_{i,j} \sim GG(\alpha,q \alpha,{\lambda}_{i}^*)$, $j=1,\cdots,n_i$, $i=1,2$.
 
 \end{itemize}
Assume that $ \lambda_1\le  \lambda^{*}_1 \le \lambda^{*}_2 \le \lambda_2$. Then, for $q \le 1 $, the following three statements are
equivalent:
\begin{itemize}
\item[(i)]$(\underbrace{\lambda_1,\cdots,\lambda_1}_{n_1},
\underbrace{\lambda_2,\cdots,\lambda_2}_{n_2})\stackrel{p}{\succeq}
(\underbrace{\lambda^{*}_1,\cdots,\lambda^{*}_1}_{n_1},\underbrace{\lambda^{*}_2,\cdots,\lambda^{*}_2}_{n_2})$;
\item[(ii)] $X_{n:n}\ge_{hr}Y_{n:n}$;
\item[(iii)]$X_{n:n}\ge_{st}Y_{n:n}$.
\end{itemize}
\et

\section*{Appendix}
Suppose $\tilde{r}$ denotes the reversed hazard rate function
corresponding to $F$. Then, from the assumption in {\bf C-1}, we
easily observe that
\bea y\tilde{r}(y)&=&\frac{\displaystyle y f(y)}{\displaystyle\int_{0}^{y}f(u)du}\nonumber\\
&=&\frac{\displaystyle
w(y)h(y)}{\displaystyle\int_{0}^{1}\,w(yu)\,h(yu)du}\nonumber\\
&=&\frac{\displaystyle h(y)}{\displaystyle
\int_{0}^{1}w(u)\,h(yu)du}.\label{e4}\eea Taking derivative from both sides
of (\ref{e4}), we obtain \bea \Big(y\tilde{r}(y)\Big)'=\frac
{\displaystyle h'(y)\int_{0}^{1}\, w(u)h(y
u)\,du-h(y)\int_{0}^{1}u\,w(u)\,h'(y u)\,du}
{\displaystyle\Big(\int_{0}^{1}\,w(u)\,h(y u)\,du\Big)^2}.\label{e5}\eea
Moreover, by using the integration by parts and the results in
 Remark 2-1, we can easily observe that
\begin{eqnarray*}
\int_{0}^{1}u\,w(u)\,h'(yu)\,du=\frac{\displaystyle 1
}{\displaystyle
y}\,\Big(h(y)-(1+w'(1))\int_{0}^{1}w(x)\,h(yu)\,du\Big).
\end{eqnarray*}
Now, from (\ref{e4}) and by substituting the above observation in
(\ref{e5}), we have \bea \Big(y\tilde{r}(y)\Big)'&=&\frac{\displaystyle
h'(y)}{\displaystyle \int_{0}^{1}\,w(u)\,h(y
u)\,du}-\frac{\displaystyle h^{2}(y)}{\displaystyle
y\Big(\int_{0}^{1}\,w(u)\,h(y u) du\Big)^2}+\frac{\displaystyle
\Big(1+w'(1)\Big)h(y)}{\displaystyle y\int_{0}^{1}\,w(u)\,h(y u)
du }\nonumber\\&=&\Big(y\tilde{r}(y)\Big)\,\frac{\displaystyle
h'(y)}{\displaystyle
h(y)}-y\tilde{r}^{2}(y)+\Big(1+w'(1)\Big)\tilde{r}(y).\label{e6}\eea 

\begin{lemma}
Consider the decomposition form in (\ref{e3}) under the assumptions
{\bf C-1} and {\bf C-2}. Suppose the following hold:
\begin{itemize}
\item[($a$)]$x\,(\tilde{r}(x)-\frac{h'(x)}{h(x)})$ is increasing in
$x\in\mathbb{R}^+$;
\item[($b$)]$\bar{F}(x)\le -\frac{w(x)\,h^2(x)}
{h'(x)}$ for every $x\in\mathbb{R}^+$.
\end{itemize}
Then,
\begin{itemize}
\item[(i)] $1 + w'(1) > 0$,
\item[(ii)]for $\alpha_i \ge 1,\,\,i=1,\cdots,n$, we have
\begin{eqnarray}
\sum_{i=1}^n\,\alpha_i y_i\tilde{r}(y_i)-y_p\Big(\tilde{r}(y_p)-
\frac{\displaystyle h'(y_p)}{\displaystyle h
(y_p)}\Big)\,\Big(1-\prod_{i=1}^{n}\,\big(F(y_i)\big)^{\alpha_i}\Big)\ge 0,\label{e7}
\end{eqnarray}
where $y_i\in\mathbb{R}^+$, $i=1,\cdots,n$, and
$y_p=\min(y_1,\cdots,y_n)$. Further, for the special case of
$n=1$, the result in (\ref{e7}) holds just under the assumption ($b$).
\end{itemize}
\label{L.2.3}
\end{lemma}
{\bf Proof}\,\,
\begin{itemize}
\item[(i)] If $h(x)$ is increasing, then the assumption (b) in not hold. So we assume that $h(x)$ is decreasing in $x$ and
 we have;
\begin{eqnarray*}
1&=&\int_{0}^{\infty} x^{w'(1)}h(x)dx\\
&=&\int_{0}^{1} x^{w'(1)}h(x)dx+\int_{1}^{\infty} x^{w'(1)}h(x)dx\\
& \ge & \int_{0}^{1} x^{w'(1)}h(x)dx\\
&\ge& \min_{0<x<1}(h(x))\int_{0}^{1} x^{w'(1)}dx\\
&=&h(1)\bigg(\frac{1}{w'(1)+1}-\lim_{x\longrightarrow 0}\frac{x^{1+w'(1)}}{1+w'(1)}\bigg)\Longrightarrow 1+w'(1)>0
\end{eqnarray*}
\item[(ii)]
From the assumptions ($b$) and ${\bf C-1}$, and
the decomposition form in (\ref{e3}) we have
\begin{eqnarray*}
\frac{F(y_i)}{\displaystyle\bar{ F}(y_i)}&\ge&
-\frac{\displaystyle h'(y_i)}{\displaystyle
w(y_i)\,h^2(y_i)}\,\int_{0}^{1}y_i\,w(y_iu)\,h(y_iu)\,du\nonumber\\
&=&-y_i\,\frac{\displaystyle h'(y_i)}{\displaystyle
h^2(y_i)}\,\int_{0}^{1}w(u)\,h(y_iu)\,du,\qquad\quad i=1,\cdots,n,
\end{eqnarray*}
which, according to  (\ref{e4}), implies
\begin{eqnarray}
\bar{F}(y_i)&\le& \frac{\displaystyle 1}{\displaystyle
1-\Big(y_i\,\frac{\displaystyle
h'(y_i)}{\displaystyle h^2(y_i)}\Big)\int_{0}^{1}w(u)\,h(y_iu)\,du}\nonumber\\
&=&\frac{\displaystyle \frac{\displaystyle h(y_i)}{\displaystyle
\int_{0}^{1}w(u)\,h(y_iu)\,du}}{\frac{\displaystyle
h(y_i)}{\displaystyle\int_{0}^{1}\,w(u)\,h(y_iu)\,du}-y_i\,\frac{\displaystyle
h'(y_i)}{\displaystyle
h(y_i)}}\nonumber\\&=&\frac{\displaystyle\tilde{r}(y_i)
}{\displaystyle\tilde{r}(y_i)-\frac{\displaystyle
h'(y_i)}{\displaystyle h(y_i)}},\qquad\quad i=1,\cdots,n.\label{e8}
\end{eqnarray}
\end{itemize}
Now, from  Proposition 2 of [balakrishnan and zhao 2013] and  (\ref{e8}), we obtain
\begin{eqnarray}
1-\prod_{i=1}^{n}\big(F(y_i)\big)^{\alpha_i}&\le&
\sum_{i=1}^{n}1-\big({F}(y_i)\big)^{\alpha_i}\nonumber\\
&\le &\sum_{i=1}^{n}{\alpha_i}\big(1-{F}(y_i)\big)\nonumber\\
&\le&\sum_{i=1}^{n}\alpha_i \frac{\displaystyle\tilde{r}(y_i)
}{\displaystyle\tilde{r}(y_i)-\frac{\displaystyle
h'(y_i)}{\displaystyle h(y_i)}}.\label{e9}
\end{eqnarray}
Note that, for the special case of $n=1$, (\ref{e9}) coincides with
(\ref{e8}) which holds just under the assumption ($b$). Now, let us
assume that $n>1$. For $i=1,\cdots,n$, set
$$\gamma_i=\frac{\displaystyle y_i\,\tilde{r}(y_i)}{\displaystyle h(y_i)},
\qquad\beta_i=\frac{\displaystyle\frac{\displaystyle
\tilde{r}(y_i)}{\displaystyle h(y_i)}}{\displaystyle
\tilde{r}(y_i)-\frac{\displaystyle h'(y_i)}{\displaystyle
h(y_i)}},\quad\delta_i=\alpha_i h(y_i).$$ So, from Proposition 1 of [balakrishnan and zhao 2013], it follows
that \bea \frac{\displaystyle\sum_{i=1}^{n}\alpha_i
y_i\tilde{r}(y_i)}{\displaystyle\sum_{i=1}^n\alpha_i\frac{\displaystyle\tilde{r}(y_i)}
{\displaystyle\tilde{r}(y_i)-\frac{\displaystyle
h'(y_i)}{\displaystyle h(y_i)}}}&\ge&\min_{1\le i\le
n}\,\Bigg\{\frac{\displaystyle\gamma_i}{\displaystyle\beta_i}\Bigg\}\nonumber\\
&=&y_p\Big(\tilde{r}(y_p)-\frac{\displaystyle
h'(y_p)}{\displaystyle h(y_p)}\Big),\label{e10}\eea where the last equality
is obtained from the assumption ($a$). Now, based on (\ref{e9}) and
(\ref{e10}), we get
\begin{eqnarray*}
\frac{\displaystyle\sum_{i=1}^n\alpha_i y_i\tilde{r}(y_i)}{\displaystyle
1-\prod_{i=1}^{n}\big(F(y_i)\big)^{\alpha_i }}&\ge&\frac{\displaystyle\sum_{i=1}^n \alpha_i 
y_i\tilde{r}(y_i)}{\displaystyle\sum_{i=1}^{n}\alpha_i \frac{\displaystyle\tilde{r}(y_i)
}{\displaystyle\tilde{r}(y_i)-\frac{\displaystyle
h'(y_i)}{\displaystyle h(y_i)}}}\nonumber\\&\ge&
y_p\Big(\tilde{r}(y_p)-\frac{\displaystyle h'(y_p)}{\displaystyle
h(y_p)}\Big),
\end{eqnarray*}
which completes the proof of the lemma.$\qquad\square$

\begin{lemma}
Suppose $X_1,\cdots,X_n$ are independent non-negative random
variables with $X_i\sim F(\lambda x)$, $i=1,\cdots,p$, and
$X_j\sim F(\lambda^*x)$, $j=p+1,\cdots,n,$ where
$\lambda\le\lambda^*$. Further, suppose the following assumptions
hold:
\begin{itemize}
\item[($a$)]$x\,(\tilde{r}(x)-\frac{h'(x)}{h(x)})$ is increasing in
$x\in\mathbb{R}^+$;
\item[($b$)]$\bar F(x)\le -\frac{
w(x)\,h^2(x)}{h'(x)}$ for every $x\in\mathbb{R}^+$.
\end{itemize}
Then, the hazard rate function of $X_{n:n}$ is increasing in
$\lambda\in(0,\lambda^*].$
\label{L.2.4}
\end{lemma}

{\bf Proof}\,\,The hazard rate function of $X_{n:n}$ can be
expressed as \ben r_{X_{n:n}}(x)&=&\Bigg(\frac{F^{p}(\lambda
x)F^{n-p}(\lambda^*x)} {1-F^{p}(\lambda
x)F^{n-p}(\lambda^*x)}\Bigg)
\Big(p\lambda\tilde{r}(\lambda x)+(n-p)\lambda^*\tilde{r}(\lambda^*x)\Big)\\
&=&\frac{\displaystyle 1}{\displaystyle x}\varphi(\lambda
x,\lambda^*x),\qquad x>0, \een where the continuously
differentiable function
$\phi:\mathbb{R^+}^2\rightarrow\mathbb{R^+}$ is defined as
\bea\varphi(y,y^*)= \Bigg(\frac{F^{p}(y)F^{n-p}(y^*)}
{1-F^{p}(y)F^{n-p}(y^*)}\Bigg) \Big(p
y\tilde{r}(y)+(n-p)y^*\tilde{r}(y^*)\Big),\quad y\le y^*.\label{e11}\eea  In
order to prove the required result, it suffices to show that
$\varphi(y,y^*)$ is increasing in $y\in(0,y^*]$.  The partial
derivative of $\varphi(y,y^*)$ with respect to $y$ is \ben
\frac{\displaystyle\partial\varphi(y,y^*)}{\displaystyle\partial
y}=p\Bigg(\frac{F^{p}(y)F^{n-p}(y^*)}
{\Big(1-F^{p}(y)F^{n-p}(y^*)\Big)^2}\Bigg)\Bigg(\tilde{r}(y)\Big(p
y_1\tilde{r}(y)+(n-p)y^*\tilde{r}(y^*)\Big)+\Big(y\tilde{r}(y)\Big)'
\,\Big(1-F^{p}(y)F^{n-p}(y^*)\Big)\Bigg),\een which, based on
(\ref{e6}), becomes \bea
\frac{\displaystyle\partial\varphi(y,y^*)}{\displaystyle\partial
y}&=&p\Bigg(\frac{F^{p}(y)F^{n-p}(y^*)}
{\Big(1-F^{p}(y)F^{n-p}(y^*)\Big)^2}\Bigg)\Bigg\{\Big(1+w'(1)\Big)\Big(1-F^{p}(y)F^{n-p}(y^*)\Big)
\tilde{r}(y)\nonumber\\
&&\quad+\tilde{r}(y)\,\Bigg(p
y\tilde{r}(y)+(n-p)y^*\tilde{r}(y^*)-y\Big(\tilde{r}(y)-\frac{\displaystyle
h'(y)}{\displaystyle
h(y)}\Big)\Big(1-F^{p}(y)F^{n-p}(y^*)\Big)\Bigg)\Bigg\}.\label{e12}\eea Now,
from the part (i) of Lemma \ref{L.2.3}, we  observe that the first term on the
right hand side of (\ref{e12}) is non-negative. On the other hand,
since $y\le y^*$, then from the assumption $(a)$ and $(b)$ and
Lemma \ref{L.2.3}, it follows that the second term on the right hand side
of (\ref{e12}) is also non-negative. This completes the proof of the
lemma. $\qquad\square$\\

\begin{lemma}
Suppose $X_1,\cdots,X_n$ are independent non-negative random
variables with common distribution function $\big(F(\lambda x)\big)^{\alpha}$.
Further, suppose the following hold:\\

$\bar F(x)\le -\frac{
w(x)\,h^2(x)}{h'(x)}$ for every $x\in\mathbb{R}^+$.\\

Then, for $\alpha \ge 1$ the hazard rate function of $X_{n:n}$ is increasing in
$\lambda\in\mathbb{R}^+$.
\label{L.2.5}
\end{lemma}

{\bf Proof}\,\,Suppose $Y_1,\cdots,Y_n$ are independent
non-negative random variables with common distribution function
$\big(F(\mu x)\big)^{\alpha_i}$ where $\lambda\le\mu$. In order to prove the required
result, it suffices to show that $X_{n:n}\ge_{hr}Y_{n:n}$. To
this end, in view of Theorem 1.B.36 of Shaked and Shanthikumar
(2007), it is enough to show that $X_i \ge_{hr}Y_i$ for
$i=1,\cdots,n$. Setting $p=n=1$ in (\ref{e11}), we observe that the
hazard rate function of $X_i$ is
$r_{X_i}(x)=\frac{1}{x}\,\varphi(\lambda x),$. Now, to derive the
result that $X_i\ge_{hr}Y_i$, we need to show that $\varphi$ is
increasing in $y\in\mathbb{R}^+$. From (\ref{e12}), we have
\begin{eqnarray}
\varphi^{\prime}(y)=\alpha\frac{\displaystyle\big( F(y)\big)^{\alpha}} {\displaystyle
\Big(1-\big( F(y)\big)^{\alpha}\Big)^2}
\Bigg\{\Big(1+w'(1)\Big)\big( F(y)\big)^{\alpha}\tilde{r}(y)+\tilde{r}(y)
\Bigg( \alpha y\tilde{r}(y)-y \Big(\tilde{r}(y)-y\frac{\displaystyle
h'(y)}{h(y)}\Big)\Big(1-\big( F(y)\big)^{\alpha}\Big)\Bigg)\Bigg\}.\label{e13}
\end{eqnarray}
according to Lemma \ref{L.2.3}, $\varphi^{\prime}$ is non-negative and the
proof is thus completed. $\qquad\square$\\
The following lemma plays a vital role in the sequel.

\begin{lemma} For $n>1$, let $\phi:\mathbb{R}_{+}^{2n}\longrightarrow \mathbb{R}_{+}$ be a continuously differentiable mapping. If the following assumptions hold:
\begin{itemize}
\item[(a)]$(\lambda_i-\lambda_j)\big(\frac{\lambda_i}{\alpha_i}\frac{\partial \phi(\alpha,\lambda)}{\partial \lambda_i}-\frac{\lambda_j}{\alpha_j}\frac{\partial \phi(\alpha,\lambda)}{\partial \lambda_j}\big)\ge 0$,
\item[(b)]$\frac{1}{\alpha_i}\frac{\partial \phi(\alpha,\lambda)}{\partial \lambda_i}|_{\lambda_i=z}=\frac{1}{\alpha_j}\frac{\partial \phi(\alpha,\lambda)}{\partial \lambda_j}|_{\lambda_j=z}$,
\end{itemize} 
Then, for any $\bold{x}=(x_1,\cdots,x_n)\in (0,\infty)^n$, the following inequality holds:
\begin{eqnarray}
\phi(\alpha,\bold{x})\ge \phi(\alpha,\bold{\tilde{x}}).
\label{e-new1}
\end{eqnarray}
Where, $\alpha=(\alpha_1,\cdots,\alpha_n)$, $\bold{\tilde{x}}=(\tilde{x},\cdots,\tilde{x})$ and $\tilde{x}=\prod_{i=1}^n x_i^{\frac{\alpha_i}{\sum_{i=1}^n\alpha_i}}$. 
\label{L.3.1}
\end{lemma}
{\bf Proof}
For a fixed vector $\bold{x}=(x_1,\cdots,x_n)\in (0,\infty)^n$ let consider $a = \min {x_i}$, $b = \max {x_i}$, $\tilde{\bold{x}}=(\tilde{x},\cdots,\tilde{x})$. Inequality (\ref{e-new1}) is an equality for $a = b$.\\
 Let us assume $a < b$. Then $\tilde{x} \in (a, b)$. We consider the compact subset $K$ of $(0,\infty)^n$:
 $$K=\big\{\bold{t}=(t_1,\cdots,t_n)\in [a,b]^n|\prod_{i=1}^n t_{i}^{\frac{\alpha_i}{\sum_{i=1}^n \alpha_i}}=\tilde{x}\big\}$$
 Clearly, $\bold{x}$ and $\tilde{\bold{x}}$ belong to K. From Weierstrass's theorem it follows that the continuous mapping $\phi$ reaches
an absolute minimum on the compact $K$ on some point $\bold{u} = (u_1, ..., u_n) \in K$.\\
Now let us assume $ \bold{u}\neq\tilde{\bold{x}}$. In this case, there exists $p,q \in \{1,2,\cdots, n\}$ such that $a\le u_p=\min {u_i}<\max {u_i}\le b$. Now we have,
$$u_{p}^{{\alpha_p}}u_{q}^{{\alpha_q}}\prod_{i\neq p,q}^n u_{i}^{{\alpha_i}}=m$$
Where $m=\tilde{x}^{\sum_{i=1}^n \alpha_i}$. The equation $z^{{\alpha_p}+{\alpha_q}}\prod_{i\neq p,q}^n u_{i}^{{\alpha_i}}=m$ has a positive solution in $z$ which is denoted
$z_1$. Clearly,
\begin{eqnarray}
z_1^{{\alpha_p}+{\alpha_q}}=u_{p}^{{\alpha_p}}u_{q}^{{\alpha_q}},\qquad\qquad u_p<z_1<u_q
\label{e-new2}
\end{eqnarray}
For $t \in [u_p, z1)$, let us consider the function $g(t)$ defined on $[u_p, z1)$ by the
relation:
$$g^{{\alpha_q}}(t)=\frac{u_{p}^{{\alpha_p}}u_{q}^{{\alpha_q}}}{t^{{\alpha_p}}}\in (z_1,u_q],\qquad \forall t \in [u_p, z1)$$
The continuously differentiable decreasing function g has the following derivative:
\begin{eqnarray}
g'(t)=-\frac{\alpha_p u_{p}^{{\alpha_p}}u_{q}^{{\alpha_q}} }{\alpha_q g^{\alpha_q-1}(t)\,t^{\alpha_p+1}}
\label{e-new3}
\end{eqnarray}
 Let $\bold{u}(t)=(u_1(t), u_2(t), \cdots, u_n(t))$, where $u_p(t)=t,\,\,u_q(t)=g(t)$ and  $u_i(t)=u_i$ when $i\neq p,q$. Now let us consider the continuously differentiable function $\psi:[u_p,z_1)\longrightarrow \mathbb{R}$, $\psi(t)=\phi(\alpha,\bold{u}(t))$.
  The derivative of $\psi(t)$ with respect to $t$, is
  \begin{eqnarray*}
  \psi'(t)&=&\frac{\partial \phi(\alpha,\bold{u}(t))}{\partial u_p(t)}\frac{\partial u_p(t)}{\partial t}+\frac{\partial \phi(\alpha,\bold{u}(t))}{\partial u_q(t)}\frac{\partial u_q(t)}{\partial t}\\
  &=&\frac{\partial \phi(\alpha,\bold{u}(t))}{\partial u_p(t)}+\frac{\partial \phi(\alpha,\bold{u}(t))}{\partial u_q(t)}g'(t)\\
  \end{eqnarray*}
  Let $t=u_p$, then $g(u_p)=u_q$. So
  \begin{eqnarray*}
  \psi'(u_p)
  &=&\frac{\partial \phi(\alpha,\bold{u}(t))}{\partial u_p(t)}|_{t=u_p}-\frac{\alpha_p u_q}{\alpha_q u_p}\frac{\partial \phi(\alpha,\bold{u}(t))}{\partial u_q(t)}|_{t=u_p}\\
  &=&\frac{\alpha_p}{u_p}\big[\frac{u_p}{\alpha_p}\frac{\partial \phi(\alpha,\bold{u}(t))}{\partial u_p(t)}|_{t=u_p}-\frac{u_q}{\alpha_q}\frac{\partial \phi(\alpha,\bold{u}(t))}{\partial u_q(t)}|_{t=u_p}\big]\\
  \end{eqnarray*}
  But, from assumptions (a) and (b), it follows that $\psi'(z)=0$ and for $u_p<z$, $\psi'(u_p)<0$.
  Hence, there exists $ \varepsilon > 0$ such that $u_p + \varepsilon < z$ and
$\psi'(t)<0$, $\forall \in [u_p, u_p + \varepsilon)$. Therefore, $\phi(\bold{u}(t)) < \phi(\bold{u}(u_p)) = \phi(\bold{u})$, for any $t \in  (u_p, u_p + \varepsilon)$. This gives the
contradiction. Then the unique minimum point of $\phi$ on K is $\bold{u}$ and the relation (\ref{e-new1}) follows.
$\qquad\square$\\
{\bf Proof of lemma \ref{L_new}}\\
{\bf (i)} $\Rightarrow$ {\bf (ii)}.  The hazard rate function of
$X_{n:n}$ can be rewritten as
\begin{eqnarray*}
r_{X_{n:n}}(x)&=&\Bigg(\frac{\displaystyle\prod_{i=1}^{n}
\big(F(\lambda_i x)\big)^{\alpha_i}}{\displaystyle 1-\prod_{i=1}^{n} \big(F(\lambda_i x)\big)^{\alpha_i}}\Bigg)\,
\sum_{i=1}^{n}\alpha_i\lambda_i\tilde{r}(\lambda_ix)\\
&=&\frac{\displaystyle 1}{\displaystyle
x}\,\psi(\lambda_1x,\cdots,\lambda_nx),\qquad x>0,
\end{eqnarray*}
where $\tilde{r}(x)={f(x)}/{F(x)}$ and the symmetric and  continuously differentiable function
$\psi:\,{\mathbb{R}^+}^n\rightarrow\,{\mathbb{R}^+} $ is defined
as
\begin{eqnarray}
\psi(y_1,\cdots,y_n)=\Bigg(\frac{\displaystyle\prod_{i=1}^{n}
\big(F(y_i )\big)^{\alpha_i}}{\displaystyle 1-\prod_{i=1}^{n} \big(F(y_i)\big)^{\alpha_i}}\Bigg)\,
\sum_{i=1}^{n}\alpha_i y_i\tilde{r}(y_i).\label{e17}
\end{eqnarray}
Similarly, the hazard rate function of $Y_{n:n}$ can be expressed
as
\begin{eqnarray*}
r_{Y_{n:n}}(x)=\frac{\displaystyle 1}{\displaystyle
x}\,\psi(\lambda x,\cdots,\lambda x),\qquad x>0.
\end{eqnarray*}
First,  we show that
$\psi(y_1,\cdots,y_n)\le\psi(\underbrace{\tilde{y},\cdots,\tilde{y}}_n)$,
where $\tilde{\lambda}=\prod_{i=1}^n \lambda_i^{\frac{\alpha_i}{\sum_{i=1}^n \alpha_i}}$. To this end, we will
utilize Lemma \ref{L.3.1} Set $y_{p}=\min(y_1,\cdots,y_n)$ and
$y_{q}=\max(y_1,\cdots,y_n)$. Then, we have
\begin{eqnarray*}
\frac{\displaystyle\partial\psi({\bf y})}{\displaystyle\partial
y_{p}}=\alpha_p\Bigg(\frac{\displaystyle\prod_{i=1}^{n}
\big(F(y_i)\big)^{\alpha_i}}{\displaystyle\Big(1-\prod_{i=1}^{n}
\big(F(y_i)\big)^{\alpha_i}\Big)^2}\Bigg)\Bigg(\tilde{r}(y_p)
\sum_{i=1}^{n}\alpha_i y_i\tilde{r}(y_i)+\Big(y_p\tilde{r}(y_p)\Big)'\Big(1-\prod_{i=1}^{n}
\big(F(y_i)\big)^{\alpha_i}\Big)\Bigg),
\end{eqnarray*}
which, according to (\ref{e6}), converts to
\begin{eqnarray*}
\frac{\displaystyle\partial\psi({\bf y})}{\displaystyle\partial
y_{p}}&=&\alpha_p \Bigg(\frac{\displaystyle\prod_{i=1}^{n}
\big(F(y_i)\big)^{\alpha_i}}{\displaystyle\Big(1-\prod_{i=1}^{n}\big(F(y_i)\big)^{\alpha_i}\Big)^2}\Bigg)
\Bigg\{\Big(1+w'(1)\Big)\Big(1-\prod_{i=1}^{n}
\big(F(y_i)\big)^{\alpha_i}\Big)\tilde{r}(y_p)\nonumber\\&&\qquad\qquad\qquad\qquad\quad+\tilde{r}(y_p)
\Bigg(\sum_{i=1}^{n}\,\alpha_i y_i\tilde{r}(y_i)-y_p
\Big(\tilde{r}(y_p)-y_p\frac{\displaystyle
h'(y_p)}{h(y_p)}\Big)\Big(1-\prod_{i=1}^{n}\big(F(y_i)\big)^{\alpha_i}\Big)\Bigg)\Bigg\}.
\end{eqnarray*}
Since the function $\psi$ is symmetric, then each partial
derivative of it has the same structure. Thus, we have
\begin{eqnarray*}
\frac{\displaystyle\partial\psi({\bf y})}{\displaystyle\partial
y_{p}}&=&\alpha_q \Bigg(\frac{\displaystyle\prod_{i=1}^{n}
\big(F(y_i)\big)^{\alpha_i}}{\displaystyle\Big(1-\prod_{i=1}^{n}\big(F(y_i)\big)^{\alpha_i}\Big)^2}\Bigg)
\Bigg\{\Big(1+w'(1)\Big)\Big(1-\prod_{i=1}^{n}
\big(F(y_i)\big)^{\alpha_i}\Big)\tilde{r}(y_q)\nonumber\\&&\qquad\qquad\qquad\qquad\quad+\tilde{r}(y_q)
\Bigg(\sum_{i=1}^{n}\,\alpha_i y_i\tilde{r}(y_i)-y_q
\Big(\tilde{r}(y_q)-y_q\frac{\displaystyle
h'(y_q)}{h(y_q)}\Big)\Big(1-\prod_{i=1}^{n}\big(F(y_i)\big)^{\alpha_i}\Big)\Bigg)\Bigg\}.
\end{eqnarray*}
Now, based on the above derivatives, we get
\begin{eqnarray}
\frac{y_p}{\alpha_p}\frac{\displaystyle
\partial\psi(\mbox{\boldmath$y$})}{\partial y_{p}}-\frac{y_q}{\alpha_q}\frac{\displaystyle
\partial\psi(\mbox{\boldmath$y$})}{\partial y_{q}}
&\stackrel{sgn}{=}&\Big(1+\omega^{\prime}(1)\Big)\Bigg(y_p\tilde{r}(y_p)-y_q\tilde{r}(y_q)\Bigg)
\Big(1-\prod_{i=1}^{n}\big(F(y_i)\big)^{\alpha_i}\Big)\nonumber\\
&&\,\,\,+y_p\tilde{r}(y_p)
\Bigg(\sum_{i=1}^n\,\alpha_i y_i\tilde{r}(y_i)-y_p\Big(\tilde{r}(y_p)-\frac{h'(y_p)}{h(y_p)}\Big)
\Big(1-\prod_{i=1}^{n}\big(F(y_i)\big)^{\alpha_i}\Big)\Bigg)
\nonumber\\&&\qquad-y_q\tilde{r}(y_q)
\Bigg(\sum_{i=1}^n\,\alpha_i y_i\tilde{r}(y_i)-y_q\Big(\tilde{r}(y_q)-\frac{h'(y_q)}{h(y_q)}\Big)
\Big(1-\prod_{i=1}^{n}\big(F(y_i)\big)^{\alpha_i}\Big)\Bigg) \nonumber\\&=&A_1+A_2,\label{e18}
\end{eqnarray}
where $a\stackrel{sgn}{=}b$ means that $a$ and $b$ have the same
sign and
$$A_1=\Big(1+\omega^{\prime}(1)\Big)\Bigg(y_p\tilde{r}(y_p)-y_q\tilde{r}(y_q)\Bigg)
\Big(1-\prod_{i=1}^{n}\big(F(y_i)\big)^{\alpha_i}\Big)$$ and
$$ A_2=y_p\tilde{r}(y_p)
\Bigg(\sum_{i=1}^n\,\alpha_i y_i\tilde{r}(y_i)-y_p\Big(\tilde{r}(y_p)-\frac{h'(y_p)}{h(y_p)}\Big)
\Big(1-\prod_{i=1}^{n}\big(F(y_i)\big)^{\alpha_i}\Big)\Bigg) -y_q\tilde{r}(y_q)
\Bigg(\sum_{i=1}^n\,\alpha_i y_i\tilde{r}(y_i)-y_q\Big(\tilde{r}(y_q)-\frac{h'(y_q)}{h(y_q)}\Big)
\Big(1-\prod_{i=1}^{n}\big(F(y_i)\big)^{\alpha_i}\Big)\Bigg).$$ From the part (i) of Lemma \ref{L.2.3} and assumptions
$(b)$, it readily follows that $A_1\ge 0$. On the other
hand, for $\alpha_i\ge 1,\,\,i=1,\cdots,n$ we have
\begin{eqnarray*}
A_2&\ge&y_q\tilde{r}(y_q)\Big(1-\prod_{i=1}^{n}\big(F(y_i)\big)^{\alpha_i}\Big)
\Bigg(y_q\Big(\tilde{r}(y_q)-\frac{h'(y_q)}{h(y_q)}\Big)
-y_p\Big(\tilde{r}(y_p)-\frac{h'(y_p)}{h(y_p)}\Big)\Bigg)\nonumber\\
&\ge& 0,
\end{eqnarray*}
where the first inequality obtains from the assumptions $(a)$,
$(b)$, $(c)$ and Lemma \ref{L.2.3}, while the second inequality follows
from the assumption ($b$). Also, it is easy to see that 
\begin{eqnarray}
\frac{1}{\alpha_p}\frac{\displaystyle
\partial\psi(\mbox{\boldmath$y$})}{\partial y_{p}}|_{y_{p}=z}-\frac{1}{\alpha_q}\frac{\displaystyle
\partial\psi(\mbox{\boldmath$y$})}{\partial y_{q}}|_{y_{p}=z}=0 \label{e-new1-1}
\end{eqnarray}
 Therefore, from Lemma \ref{L.3.1}, (\ref{e18}) and (\ref{e-new1-1}),
it follows that $\psi(\lambda_1x,\cdots,\lambda_nx)
\le\psi(\tilde{\lambda}x,\cdots,\tilde{\lambda}x)$. On the other
hand, from the part (i) of Lemma \ref{L.2.3} and assumption ($c$), and Lemma 2.5 it
follows that the hazard rate function of $Y_{n:n}$ is increasing
in $\lambda\in\mathbb{R}^+$, that is,
$\psi(\tilde{\lambda}x,\cdots,\tilde{\lambda}x) \le\psi(\lambda
x,\cdots,\lambda x)$ for $\lambda\ge\tilde{\lambda}$. By
combining these observations,  (ii) follows.

{\bf (ii)} $\Rightarrow$ {\bf (iii)}. It is clear since  the
hazard rate order implies the usual stochastic order.

{\bf (iii)} $\Rightarrow$ {\bf (i)}. From the statement (iii), we
have $\frac{F_{X_{n:n}}(x)}{F_{Y_{n:n}}(x)}\le 1$ for every
$x>0$. Thus,
\begin{eqnarray}
\prod_{i=1}^{n}\left(lim_{x\rightarrow
0}\frac{F(\lambda_ix)}{F(\lambda x)}\right)^{\alpha_i}&=&lim_{x\rightarrow
0}\,\prod_{i=1}^n\,\, \left(\frac{F(\lambda_i x)}{F(\lambda
x)}\right )^{\alpha_i}\nonumber\\&\le& 1.\label{e19}
\end{eqnarray}
On the other hand, according to the L'Hopital's rule, and Remark 2.1, and assumption {\bf C-2}, we obtain
\begin{eqnarray}
lim_{x\rightarrow 0}\frac{F(\lambda_i x)}{F(\lambda
x)}&=&lim_{x\rightarrow 0}\frac{\lambda_i\,f(\lambda_i x)}
{\lambda\,f(\lambda x)}\nonumber\\
&=&\Big(\frac{\lambda_i}{\lambda}\Big)^{1+w'(1)}\,\,
lim_{x\rightarrow 0}\frac{h(\lambda_ix)}{h(\lambda x)}\nonumber\\
&=&(\frac{\lambda_i}{\lambda})^{(1+w'(1))}.\label{e20}
\end{eqnarray}

Hence, by substituting (\ref{e20}) into (\ref{e19}) and using the part (i) of Lemma \ref{L.2.3}, we have $\lambda\ge\tilde{\lambda}$ which completes the
proof of the theorem. $\qquad\square$\\
 
 The following
known result provides an approach for testing whether a vector
valued function is Schur-convex (Schur-concave) or not. \bl
(Marshall et al. 2011, p. 84)\,Suppose $I\subset\mathbb{R}$ is an
open interval and suppose $\phi\,:\,I^n\rightarrow\mathbb{R}$ is
continuously differentiable. Necessary and sufficient conditions
for $\phi$ to be Schur-convex  (Schur-concave) on $I^n$ are the
symmetric property of $\phi$ on $I^n$, and  \ben
(x_i-x_j)\left(\frac{\partial \phi}{\partial x_i}({\bf
x})-\frac{\partial \phi}{\partial x_j}({\bf x})\right)\,\ge (\le)
0,\qquad \,for \,\,all\,\,i\not=j\,\,and\,\, {\bf
x}=(x_1,\cdots,x_n)\in I^n, \een where $\frac{\partial
\phi}{\partial x_i}({\bf x})$ denotes the partial derivative of
$\phi$ with respect to its $i$-th argument.
\label{marshall}
\el

{\bf Proof of lemma \ref{l-new2}}\,\,Set $a=\log\lambda$ and $a^*=\log\lambda^*$. Then,
the hazard rate function of $X_{n:n}$, in terms of $a$ and $a^*$,
is given by
\begin{eqnarray*}
r_{X_{n:n}}(x;a,a^*)=\frac{\displaystyle 1}{\displaystyle
x}\,\phi(\underbrace{xe^{a},\cdots,xe^{a}}_{p},
\underbrace{xe^{a^*},\cdots,xe^{a^*}}_{n-p}),
\end{eqnarray*}
where
\ben\phi(\underbrace{y,\cdots,y}_{p},\underbrace{y^*,\cdots,y^*}_{n-p})=
\Bigg(\frac{F^{p}(y)F^{n-p}(y^*)} {1-F^{p}(y)F^{n-p}(y^*)}\Bigg)
\Big(p y\tilde{r}(y)+(n-p)y^*\tilde{r}(y^*)\Big).\een  For
simplicity in presentation, we denote
$\phi(\underbrace{y,\cdots,y}_{p},
\underbrace{y^*,\cdots,y^*}_{n-p})$ by $\phi(y,y^*)$. In order to
prove the desired result, in the light of Lemma \ref{marshall}, it suffices
to show that for fixed $x>0$, \bea
(a-a^*)\Bigg(\frac{\displaystyle\partial\phi(xe^{a},xe^{a^*})}
{\displaystyle\partial a}-
\frac{\displaystyle\partial\psi(xe^{a},xe^{a^*})}
{\displaystyle\partial a^*}\Bigg)\le 0.\label{e21}\eea Note that
\ben\frac{\displaystyle\partial\phi(xe^{a},xe^{a^*})}{\partial
a}&=&\frac{\displaystyle\partial
\Big(xe^{a}\Big)}{\displaystyle\partial a}\,\,
{\displaystyle\partial\phi(xe^{a},xe^{a^*})}
{\displaystyle\partial\Big(xe^{a}\Big)}\\&=&xe^{a}\,\,
\frac{\displaystyle\partial\phi(xe^{a},xe^{a^*})}
{\displaystyle\partial\Big(xe^{a}\Big)}.\een Similarly,
\ben\frac{\displaystyle\partial\phi(xe^{a},xe^{a^*})}{\partial
a^*}=xe^{a^*}\,\,
\frac{\displaystyle\partial\phi(xe^{a},xe^{a^*})}
{\displaystyle\partial\Big(xe^{a^*}\Big)}.\een Now, using these
observations, we get \bea
\frac{\displaystyle\partial\phi(xe^{a},xe^{a^*})}
{\displaystyle\partial a}-
\frac{\displaystyle\partial\phi(xe^{a},xe^{a^*})}
{\displaystyle\partial a^*}=& xe^{a}\,\,
\frac{\displaystyle\partial\phi(xe^{a},xe^{a^*})}
{\displaystyle\partial\Big(xe^{a}\Big)}-xe^{a^*}\,\,
\frac{\displaystyle\partial\phi(xe^{a},xe^{a^*})}
{\displaystyle\partial\Big(xe^{a^*}\Big)}.\label{e22}\eea So, from (\ref{e22}), we
easily observe that (\ref{e21}) satisfies if we could show that  \bea
(\ln y-\ln y^*)\Bigg(y\,\frac{\displaystyle\partial\psi(y,y^*)}
{\displaystyle\partial
y}-y^*\,\frac{\displaystyle\partial\psi(y,y^*)}
{\displaystyle\partial y^*}\Bigg)\le 0.\label{e23}\eea To this end, we have
\begin{eqnarray}
y\frac{\displaystyle\partial\phi(y,y^*)}{\displaystyle\partial
y}-y^*\frac{\displaystyle\partial\psi(y,y^*)}{\displaystyle\partial
y^*}&\stackrel{sgn}{=}&(1+\omega^{\prime}(1))\Big(y\tilde{r}(y)-y^*\tilde{r}(y^*)\Big)
\Big(1-F^{p}(y)F^{n-p}(y^*)\Big)\nonumber\\
\nonumber\\&&\,\,\,+y\tilde{r}(y)\Bigg(\Big(p\,y\tilde{r}(y)+(n-p)\,y^*\tilde{r}(y^*)\Big)-y\Big
(\tilde{r}(y)-\frac{h'(y)}{h(y)}\Big)\Big(1-F^{p}(y)F^{n-p}(y^*)\Big)\Bigg)\nonumber\\
&&\quad-y^*\tilde{r}(y^*)\Bigg(\Big(p\,y\tilde{r}(y)+(n-p)\,y^*\tilde{r}(y^*)\Big)-y^*\Big
(\tilde{r}(y^*)-\frac{h'(y^*)}{h(y^*)}\Big)\Big(1-F^{p}(y)F^{n-p}(y^*)\Big)\Bigg)\nonumber\\
\nonumber\\&=&B_1+B_2,\label{e24}
\end{eqnarray}
where \ben
B_1=(1+\omega^{\prime}(1))\Big(1-F^{p}(y)F^{n-p}(y^*)\Big)\Big(y\tilde{r}(y)-y^*\tilde{r}(y^*)\Big),\een
\ben
B_2&=&y\tilde{r}(y)\Bigg(\Big(p\,y\tilde{r}(y)+(n-p)\,y^*\tilde{r}(y^*)\Big)-y\Big
(\tilde{r}(y)-\frac{h'(y)}{h(y)}\Big)\Big(1-F^{p}(y)F^{n-p}(y^*)\Big)\Bigg)\nonumber\\&&\quad
-y^*\tilde{r}(y^*)\Bigg(\Big(p\,y\tilde{r}(y)+(n-p)\,y^*\tilde{r}(y^*)\Big)-y^*\Big
(\tilde{r}(y^*)-\frac{h'(y^*)}{h(y^*)}\Big)\Big(1-F^{p}(y)F^{n-p}(y^*)\Big)\Bigg).\een
Let us first assume that $y^*\ge y$. Then, from assumptions ($b$)
and  the part (i) of Lemma \ref{L.2.3}, we can easily observe that $B_1\ge 0$. Furthermore,
according to assumptions $(a)$, $(b)$, ($c$) and Lemma \ref{L.2.3}, it
follows that
\begin{eqnarray*}
B_2&\ge&y^*\tilde{r}(y^*)\Big(1-F^{p}(y)F^{n-p}(y^*)\Big)\Bigg(y^*\Big
(\tilde{r}(y^*)-\frac{h'(y^*)}{h(y^*)}\Big)-y\Big
(\tilde{r}(y)-\frac{h'(y)}{h(y)}\Big)\Bigg)\\
&\ge& 0.
\end{eqnarray*}
These observations result that (\ref{e24}) is non-negative. For the
case when $y>y^*$, by using an argument similar to the above, one
can easily show that $B_1\le 0$ and $B_2\le 0$, i.e, (\ref{e24}) is
non-positive. Consequently, (\ref{e23}) is satisfied and the proof is
thus completed. $\qquad\square$\\
{\bf  Proof of lemma \ref{l-new3}}\,\,Suppose (i) holds. Then,  (ii) immediately follows
from Lemma \ref{l-new2}
 for the case when
$\lambda^p{\lambda^*}^{n-p}=\mu^p{\mu^*}^{n-p}$. Now, let us
assume that $\lambda^p{\lambda^*}^{n-p}<\mu^p{\mu^*}^{n-p}$.
Setting $\lambda_0=\mu\,(\frac{\mu^*}{\lambda^*})^{(n-p)/p}$, then
it can be easily seen that \ben
\lambda<\lambda_0\le\mu\le\mu^*\le\lambda^*\qquad and\qquad
\lambda^p_0{\lambda^*}^{n-p}=\mu^p{\mu^*}^{n-p}.\een Suppose
$Z_1,\cdots,Z_n$ are independent non-negative random variables
with $Z_i\sim F(\lambda_0x)$, $i=1,\cdots,p$, and $Z_j\sim
F(\lambda^*x)$, $j=p+1,\cdots,n$. From Lemma \ref{l-new2}, we have
$Z_{n:n}\ge_{hr}Y_{n:n}$. On the other hand, from Lemma \ref{L.2.4}, it
follows that $X_{n:n}\ge_{hr}Z_{n:n}$. Now, by combining these
results, (ii) follows. Since the hazard rate order implies the
usual stochastic order, then (ii) $\Rightarrow$ (iii). The proof
of the implication (iii) $\Rightarrow$ (i) is similar to that of
Lemma \ref{L_new}, and therefore omitted here. This completes the proof
of the theorem.$\qquad\square$\\
\bl
Let the function $\psi(.):\mathbb{R}^+\longrightarrow \mathbb{R}^+$ be defined as 
\begin{eqnarray*}
	\psi(x)=q\,(q-1)+e^{x-1}\bigg((q-1)^2x+q\,x^{2}-x^{1+q}+2q(1-q)\bigg)
\end{eqnarray*}
Then, for each $q\in (0,1]$, $\psi(x)\ge 0$. 
\label{EWG-2}
\el
{\bf Proof.} If $q=1$, we can easily observe that $\psi(x)= 0$ for all $ x \in [0,\infty)$. Now, let us assume that $q<1$. 
After some simplifications, we obtain that
\begin{eqnarray}
\frac{\partial}{\partial x}\psi(x)|_{
x=1}
 &=&0\nonumber\\
\frac{\partial^2}{\partial x^2}\psi(x)|_{
	x=1} &=&0 \nonumber
\end{eqnarray}
And
\begin{eqnarray*}
	\frac{\partial^3}{\partial x^3}\psi(x)|_{
		x=1} &=&-q^3-q^2+2q\\
	&>& 0\qquad \forall q\in(0,1).
\end{eqnarray*}
Therefore, $\psi'(x)$ has a local minimum at $x=1$ for all $q \in (0,1)$. So,  there exists some $\epsilon > 0$ such that
$\psi(x)>\psi(1)$ for all $x \in (1,1+\epsilon)$ and $q \in (0,1)$. Also, it is easy to see that
\begin{eqnarray*}
	\psi(1)=0\qquad and \qquad \lim_{x\longrightarrow \infty}\psi(x)=\infty
	.\end{eqnarray*}
Now, suppose $\min_{x\ge 1} \psi(x) <0$. Then, From the above observation, $\psi(x)$ has to cross the level $0$ at least twice. That is, there exist a point $\xi >0$ such that $\psi(\xi)=0$. Then using this in
\begin{eqnarray*}
	\frac{\partial}{\partial x}\psi(x)=e^{x-1}\bigg((q-1)^2x+q\,x^{2}-x^{1+q}+2q(1-q)\bigg)+e^{x-1}\bigg((q-1)^2+2\,q\,x-(1+q)\,x^{q}\bigg)
\end{eqnarray*}
we obtain that,
\begin{eqnarray*}
	\frac{\partial}{\partial x}\psi(x)|_{
		x=\xi}
 =e^{\xi-1}\bigg((q-1)^2+2\,q\,\xi-(1+q)\xi^{q}\bigg)-q\,(q-1).
\end{eqnarray*}
Let,
\begin{eqnarray*}
	g(x)=e^{x-1}\bigg((q-1)^2+2\,q\,x-(1+q)x^{q}\bigg)-q\,(q-1)
	.\end{eqnarray*}
Now, we need to show that $g(x) \ge 0$ for all $x\in [1,\infty)$. It is easy to see that
\begin{eqnarray*}
	g(1)=0\qquad and \qquad \lim_{x\longrightarrow \infty}g(x)=\infty.
	.\end{eqnarray*}
Suppose $\min_{x>1} g(x)<0 $. The derivative of $g(x)$ with respect to $x$, is
\begin{eqnarray*}
	g'(x)=e^{x-1}\bigg((q-1)^2+2\,q\,x-(1+q)x^{q}\bigg)+e^{x-1}\bigg(2\,q-q\,(1+q)x^{q-1}\bigg),
	.\end{eqnarray*}
Using the above relations, $g'(1)=2\,b^2 >0$. Then $g(x)$ has to cross the level $0$ at least twice, since $g(1)=0$, $g(\infty)=\infty$, and $g'(1)=2\,b^2 >0$. At any point $t$ that $g(t)=0$, we have
\begin{eqnarray*}
	e^{t-1}\bigg((q-1)^2+2\,q\,t-(1+q)t^{q}\bigg)=q\,(q-1)
\end{eqnarray*}
Using this in $g'(x)$,
\begin{eqnarray*}
	g'(x)|_{
		x=t} =e^{t-1}\bigg(2\,q-q(1+q)t^{q-1}\bigg)+q\,(q-1)
\end{eqnarray*} 
It is easy to see that $g'(t)>0$. Thus, $t$ is unique. That is a contradiction with  $g(x)$ has to cross the level $0$ at least twice. Then $\psi' (\xi)>0$ , since $g(x)>0$ for all $x \in (1, \infty)$. Thus, $\xi$ is unique. That is a contradiction with  $\psi(x)$ has to cross the level $0$ at least twice. So, the crossing point $\xi$ does not exist, from which the required result follows.
$\qquad\square$\\

{\bf Proof of lemma \ref{lemma_3.1}.}
As mentioned before, the baseline density function can be written
as the decomposition form in (\ref{e3}) with
\ben w(x)=x^{p-1}\qquad and\qquad
h(x)=\frac{\displaystyle p}{\displaystyle q}(1+x^p)^{\frac{1}{q}-1}e^{1-(1+x^p)^{\frac{1}{q}}},\qquad
x\in\mathbb{R^+},\een
\citep{khaledi2011stochastic} showed that $x\tilde{r}(x)$ is decreasing in
$x\in\mathbb{R^+}$, and so the assumption ($a$) of Lemma \ref{L_new} is
satisfied. Now, let us
check the assumption ($b$). After some
algebraic computation, we obtain 
\ben
s(x)&=&x\Big(\tilde{r}(x)-\frac{\displaystyle h'(x)}{\displaystyle
	h(x)}\Big)\\&=&\frac{p}{q}\,x^p\bigg(\frac{\displaystyle(1+x^p)^{\frac{1}{q}-1}\,e^{1-(1+x^p)^{\frac{1}{q}}}}{\displaystyle 1-e^{1-(1+x^p)^{\frac{1}{q}}}}-\frac{\displaystyle 1-q-(1+x^p)^{\frac{1}{q}}}{\displaystyle 1+x^p}\bigg)\\
&=&\frac{p}{q}h(1+x^p).\een
Where the continuously differentiable $h:[1,\infty)\longrightarrow \mathbb{R}^+$ is defined as
\ben
h(x)&=&(x-1)\bigg(\frac{\displaystyle x^{\frac{1}{q}-1}}{\displaystyle e^{x^{\frac{1}{q}}-1}-1}-\frac{\displaystyle 1-q-x^{\frac{1}{q}}}{\displaystyle x}\bigg)
.\een
Then, we have
\ben
\frac{\partial s(x)}{\partial x}&=&\frac{\partial (1+x^p)}{\partial x}\frac{\partial s(x)}{\partial (1+x^p)}\\
&\overset{sgn}{=}&\frac{\partial h(x)}{\partial x}\\
&=&\frac{\displaystyle x^{\frac{1}{q}-1}}{\displaystyle e^{x^{\frac{1}{q}}-1}-1}-\frac{\displaystyle 1-q-x^{\frac{1}{q}}}{\displaystyle x}
+(x-1)\bigg(\frac{\displaystyle x^{\frac{1}{q}-1}(\frac{1}{q}-1)}{\displaystyle x\,\big(e^{x^{\frac{1}{q}}-1}-1\big)}-\frac{\displaystyle x^{\frac{2}{q}-2}\,e^{x^{\frac{1}{q}}-1}}{\displaystyle q\, \big(e^{x^{\frac{1}{q}}-1}-1\big)^2}+\frac{\displaystyle x^{\frac{1}{q}}}{\displaystyle q\,x^2}+\frac{\displaystyle 1-q-x^{\frac{1}{q}}}{\displaystyle x^2}\bigg)\\
&=&\frac{1}{\displaystyle (e^{x^{\frac{1}{q}}-1}-1)^2\,b\,x^2}\big(-q+q^2e^{2t^{\frac{1}{q}}-2}-2\,q^2e^{x^{\frac{1}{q}}-1}-x^{\frac{q+2}{q}}e^{x^{\frac{1}{q}}-1}+q^2-x^{\frac{q+1}{q}}e^{x^{\frac{1}{q}}-1}-x^{\frac{1}{q}}q\,e^{x^{\frac{1}{q}}-1}\\
&&+x^{\frac{1}{q}}e^{x^{\frac{1}{q}}-1}+x^{\frac{q+1}{q}}e^{2x^{\frac{1}{q}}-2}+x^{\frac{1}{q}}q\,e^{2x^{\frac{1}{q}}-2}-x^{\frac{1}{q}}e^{2x^{\frac{1}{q}}-2}+2\,q\,e^{x^{\frac{1}{q}}-1}+x^{\frac{2}{q}}e^{x^{\frac{1}{q}}-1}-q\,e^{2x^{\frac{1}{q}}-2}\big)\\
&\overset{sgn}{=}&q\,(q-1)+e^{2x^{\frac{1}{q}}-2}\big(q(q-1)+x^{\frac{1+q}{q}}+x^{\frac{1}{q}}(q-1)\big)+
e^{x^{\frac{1}{q}}-1}\big(x^{\frac{2}{q}}-x^{\frac{1+q}{q}}-x^{\frac{2+q}{q}}+(1-q)x^{\frac{1}{q}}+2\,q\,(1-q)\big)\\
&\ge &q\,(q-1)+e^{x^{\frac{1}{q}}-1}\bigg(x^{\frac{1}{q}}\big(q(q-1)+x^{\frac{1+q}{q}}+x^{\frac{1}{q}}(q-1)\big)+x^{\frac{2}{q}}-x^{\frac{1+q}{q}}-x^{\frac{2+q}{q}}+(1-q)x^{\frac{1}{q}}+2\,q\,(1-q)
\bigg)\\
&= &
q\,(q-1)+e^{x^{\frac{1}{q}}-1}\bigg((q-1)^2x^{\frac{1}{q}}+q\,x^{\frac{2}{q}}-x^{\frac{1+q}{q}}+2\,q\,(1-q)
\bigg)\\
&=&\psi(x^{\frac{1}{q}})
.\een
Where the continuously differentiable $\psi:[1,\infty)\longrightarrow \mathbb{R}^+$ is defined as
\begin{eqnarray*}
	\psi(x)=q\,(q-1)+e^{x-1}\bigg((q-1)^2x)+qx^{2}-x^{1+q}+2q(1-q)\bigg),
\end{eqnarray*}
Now, using Lemma \ref{EWG-2}
\begin{eqnarray*}
	\psi(x)\ge 0 \qquad \forall x\in [1,\infty)\qquad and \qquad \forall q\in (0,1]
	,\end{eqnarray*}
Therefore, the assumption
($b$) of Lemma \ref{L_new} is satisfied for $0<q\le 1$. Moreover,
for $0<q\le 1$, we have
\ben
\bar{F}(x)&=&e^{1-(1+x^{p})^{\frac{1}{q}}}\\
&\le &-\frac{(1+x^{p})^{\frac{1}{q}}e^{1-(1+x^{p})^{\frac{1}{q}}}}{1-q-(1+x^{p})^{\frac{1}{q}}}\\
&=&-\frac{w(x)h^2(x)}{h'(x)}\qquad x\ge 0.
\een
which confirms the assumption ($c$) of Lemma \ref{L_new}.
This completes the proof of the lemma. $\qquad\square$\\
{\bf Proof of lemma \ref{Theo_3.2}.}
As mentioned before, the baseline density function can be written
as the decomposition form in (\ref{e3}) with \ben
w(x)=x^{\beta-1}\qquad and\qquad
h(x)=\frac{\displaystyle\alpha}{\displaystyle\Gamma^*(\frac
	{\displaystyle\beta}{\displaystyle\alpha})}\,e^{-x^{\alpha}}.\een
\citep{khaledi2011stochastic} showed that $x\tilde{r}(x)$ is decreasing in
$x\in\mathbb{R^+}$, and so the assumption ($a$) of  Lemma \ref{L_new} is
satisfied. Now, let us
check the assumption ($b$). After some
algebraic computation, we obtain \ben
s(x)&=&x\Big(\tilde{r}(x)-\frac{\displaystyle h'(x)}{\displaystyle
	h(x)}\Big)\\&=&\alpha\,x^{\alpha}+\frac{\displaystyle 1
}{\displaystyle\int_{0}^{1}u^{\beta-1}\,e^{(1-u^{\alpha})x^{\alpha}}du}.\een
Taking derivative from $s(x)$ with respect to $x$ gives rise to
\ben s'(x)&=&\alpha^2\,x^{\alpha-1}-\alpha\,x^{\alpha-1}
\frac{\displaystyle
	\int_{0}^{1}u^{\beta-1}(1-u^{\alpha})\,e^{(1-u^{\alpha})x^{\alpha}}du}{\displaystyle\Big(
	\int_{0}^{1}u^{\beta-1}\,e^{(1-u^{\alpha})x^{\alpha}}du\Big)^2}\\
&\ge&\alpha^2\,x^{\alpha-1}-\frac{\displaystyle\alpha\,x^{\alpha-1}}
{\displaystyle\int_{0}^{1}u^{\beta-1}\,e^{(1-u^{\alpha})x^{\alpha}}du}\\
&\ge&\alpha^2\,x^{\alpha-1}-\frac{\displaystyle\alpha\,x^{\alpha-1}}
{\displaystyle\int_{0}^{1}u^{\beta-1}du}\\&=&\alpha(\alpha-\beta)\,x^{\alpha-1}.\een
Thus, we can easily observe that $s(x)$ is increasing in
$x\in\mathbb{R^+}$ for $\alpha\ge\beta$. Therefore, the assumption
($b$) of  Lemma \ref{L_new} is satisfied for $\alpha\ge\beta$. Moreover,
for $\alpha\ge\beta$, we have
\ben\bar{F}(x)&=&\int_{x}^{\infty}\frac{\alpha}
{\Gamma(\frac{\beta}{\alpha})}\,u^{\beta-1}\,e^{-u^{\alpha}}du\\
&=&\int_{x}^{\infty}\frac{\alpha}
{\Gamma(\frac{\beta}{\alpha})}\,u^{\beta-\alpha}\,u^{\alpha-1}\,
e^{-u^{\alpha}}du\\&\le&\int_{x}^{\infty}\frac{\alpha\,x^{\beta-\alpha}}
{\Gamma(\frac{\beta}{\alpha})}\,u^{\alpha-1}\,e^{-u^{\alpha}}du\\&=&
\frac{x^{\beta-\alpha}}
{\Gamma(\frac{\beta}{\alpha})}\,e^{-x^{\alpha}}\\&=&-\frac{w(x)\,h^2(x)}{h'(x)},\qquad
x>0, \een which confirms the assumption ($c$) of  Lemma \ref{L_new}
This completes the proof of the lemma. $\qquad\square$\\

\section*{References}
\bibliographystyle{plain}
\bibliography{ref.bib}
\end{document}